\documentclass[number,seceqn,secthm,citesort,dvips]{arxbj}
\usepackage{cursive}   

\aid{0}
\volume{16}
\issue{2}
\pubyear{2010}
\firstpage{361}
\lastpage{388}
\doi{10.3150/09-BEJ224}

\makeatletter

\newcommand{\indic}{\mathbb{I}}
\newremark{rem}{Remark}[section]

\newtheorem{prop}{Proposition}[section]
\newtheorem{lem}{Lemma}[section]

\makeatother

\begin{document}
\begin{frontmatter}

\title{Dirichlet mean identities and laws of a class of subordinators}
\runtitle{Dirichlet mean identities}

\begin{aug}
\author{\fnms{Lancelot F.} \snm{James}\ead[label=e1]{lancelot@ust.hk}\corref{}}
\runauthor{L. F. James}
\address{Department of Information Systems, Business Statistics and
Operations Management,
The Hong Kong University of Science and
Technology, Clear Water Bay, Kowloon, Hong Kong. \printead{e1}}
\end{aug}

\received{\smonth{5} \syear{2008}}

%
\begin{abstract}
An interesting line of research is the investigation of the laws
of random variables known as Dirichlet means. However, there is not
much information on interrelationships between different
Dirichlet means. Here, we introduce two distributional operations,
one of which consists of multiplying a mean functional by an independent
beta random variable, the other being an operation involving an exponential
change of measure. These operations identify relationships between
different means and their densities. This allows one to use the
often considerable analytic work on obtaining results for one
Dirichlet mean to obtain results for an entire family of otherwise
seemingly unrelated Dirichlet means. Additionally, it allows one to
obtain explicit densities for the related class of random
variables that have generalized gamma convolution distributions
and the finite-dimensional distribution of their associated L\'{e}vy
processes. The importance of this latter statement is that L\'{e}vy
processes now commonly appear in a variety of applications in
probability and statistics, but there are relatively few cases
where the relevant densities have been described explicitly. We
demonstrate how the technique allows one to obtain the
finite-dimensional distribution of several interesting
subordinators which have recently appeared in the literature.
\end{abstract}

%

\end{frontmatter}

\section{Introduction}
In this work, we present two distributional operations which
identify relationships between seemingly different classes of
random variables which are representable as linear functionals of
a~Dirichlet process, otherwise known as \textit{Dirichlet means}.
Specifically, the first operation consists of multiplication of a
Dirichlet mean by an independent beta random variable and the
second operation involves an exponential change of measure to the
density of a related infinitely divisible random variable having a
generalized gamma convolution distribution~(GGC). This latter
operation is often referred to in the statistical literature as
\textit{exponential tilting} or in mathematical finance as an
\textit{Esscher transform}. We believe our results add a significant
component to the foundational work of Cifarelli and
Regazzini~\cite{CifarelliRegazzini79,CifarelliRegazzini}. In
particular, our results allow one to use the often considerable
analytic work on obtaining results for one Dirichlet mean to obtain
results for an entire family of otherwise seemingly unrelated mean
functionals. It also allows one to obtain explicit densities for
the related class of infinitely divisible random variables which
are generalized gamma convolutions and an explicit description of
the finite-dimensional distribution of their associated L\'{e}vy
processes (see Bertoin~\cite{Bertoin} for the formalities of
general L\'{e}vy processes). The importance of this latter statement
is that L\'{e}vy processes now commonly appear in a variety of
applications in probability and statistics, but there are
relatively few cases where the relevant
densities have been described explicitly. A detailed summary and
outline of our
results may be found in Section~\ref{sec:outline}. Some background
information on, and notation for, Dirichlet processes and Dirichlet
means, their connection with GGC random variables, recent
references and some motivation for our work are given in the next
section.

\subsection{Background and motivation}
Let $X$ be a non-negative random variable with cumulative
distribution function $F_{X}$. Note, furthermore, that for a measurable
set $C,$ we use the notation $F_{X}(C)$ to mean the probability
that $X$ is in $C.$ One may define a Dirichlet process random
probability measure (see Freedman~\cite{Freedman} and
Ferguson~\cite{Ferguson73,Ferguson74}), say $P_{\theta},$ on
$[0,\infty)$ with total mass parameter $\theta$ and prior
parameter $F_{X},$ via its finite-dimensional distribution as
follows: for any disjoint partition on $[0,\infty)$, say
$(C_{1},\ldots, C_{k})$, the distribution of the random vector
$(P_{\theta}(C_{1}),\ldots,P_{\theta}(C_{k}))$ is a $k$-variate
Dirichlet distribution with parameters $(\theta
F_{X}(C_{1}),\ldots, \theta F_{X}(C_{k})).$ Hence, for each $C$,
\[
P_{\theta}(C)=\int_{0}^{\infty}\indic(x\in C)P_{\theta}(\mathrm{d}x)
\]
has a beta distribution with
parameters $(\theta F_{X}(C),\theta(1-F_{X}(C))).$ Equivalently,
setting\break $\theta F_{X}(C_{i})=\theta_{i}$ for $i=1,\ldots,k,$
\[
(P_{\theta}(C_{1}),\ldots, P_{\theta}(C_{k}))
\stackrel{\mathrm{d}}=\biggl(\frac{G_{\theta_{i}}}{G_{\theta}};
i=1,\ldots,k\biggr),
\]
where $(G_{\theta_{i}})$ are independent random variables with
gamma$(\theta_{i},1)$ distributions and
$G_{\theta}=G_{\theta_{1}}+\cdots+G_{\theta_{k}}$ has a
gamma$(\theta,1)$ distribution. This means that one can define the
Dirichlet process via the normalization of an independent
increment gamma process on $[0,\infty)$, say
$\gamma_{\theta}(\cdot),$ as
\[
P_{\theta}(\cdot)=\frac{\gamma_{\theta}(\cdot)}{\gamma_{\theta
}([0,\infty))},
\]
where $\gamma_{\theta}(C_{i})\stackrel{\mathrm{d}}=G_{\theta_{i}}$, whose
almost surely finite total random mass is
$\gamma_{\theta}([0,\infty))\stackrel{\mathrm{d}}=G_{\theta}.$ A very
important aspect of this construction is the fact that
$G_{\theta}$ is independent of $P_{\theta}$ and hence of any
functional of $P_{\theta}.$ This is a natural generalization of
Lukacs' \cite{Lukacs} characterization of beta and gamma random
variables, which is fundamental to what is now referred to as
the \textit{beta--gamma algebra} (for more on this, see Chaumont and
Yor~(\cite{Chaumont}, Section 4.2); see also Emery and
Yor~\cite{EmeryYor} for some interesting relationships between
gamma processes, Dirichlet processes and Brownian bridges).
Hereafter, for a random probability measure $P$ on $[0, \infty),$
we write
\[
P\sim\Pi_{\theta,F_{X}},
\]
to indicate that $P$ is a
Dirichlet process with parameters $(\theta,F_{X}).$

These simple representations and other nice features of the
Dirichlet process have, since the important work of
Ferguson~\cite{Ferguson73,Ferguson74}, contributed greatly to the
relevance and practical utility of the field of Bayesian non- and
semi-parametric statistics. Naturally, owing to the ubiquity of
the gamma and beta random variables, the Dirichlet process also
arises in other areas. One of the more interesting and, we believe,
quite important topics related to the Dirichlet process is the
study of the laws of random variables called \textit{Dirichlet mean
functionals}, or simply Dirichlet means, which we denote as
\[
M_{\theta}(F_{X})\stackrel{\mathrm{d}}=\int_{0}^{\infty}xP_{\theta}(\mathrm{d}x),
\]
as initiated in the works of Cifarelli and
Regazzini~\cite{CifarelliRegazzini79,CifarelliRegazzini}.
In~\cite{CifarelliRegazzini}, the authors obtained an important
identity for the Cauchy--Stieltjes transform of order $\theta.$
This identity is often referred to as the \textit{Markov--Krein
identity}, as
can be seen in, for example, Diaconis and
Kemperman~\cite{Diaconis}, Kerov~\cite{Kerov} and Vershik, Yor and
Tsilevich~\cite{Vershik}, where these authors highlight its
importance to, for instance, the study of the Markov moment
problem, continued fraction theory and exponential representation
of analytic functions. This identity is later called the
\textit{Cifarelli--Regazzini identity} in~\cite{James2005}. Cifarelli and
Regazzini~\cite{CifarelliRegazzini}, owing to their primary
interest, used this identity to then obtain explicit density and
cdf formulae for $M_{\theta}(F_{X}).$ The density formulae may be
seen as Abel-type transforms and hence do not always have simple
forms, although we stress that they are still useful for some
analytic calculations. The general exception is the case
$\theta=1$, which has a nice form. Some examples of works that have
proceeded along these lines are Cifarelli and
Melilli~\cite{CifarelliMelilli}, Regazzini, Guglielmi and di
Nunno~\cite{Regazzini2002}, Regazzini, Lijoi and
Pr\"{u}nster \cite{Regazzini2003}, Hjort and Ongaro~\cite{Hjort},
Lijoi and Regazzini~\cite{Lijoi}, and Epifani, Guglielmi and
Melilli~\cite{Epifani2004,Epifani2006}. Moreover, the recent
works of~Bertoin \textit{et al.}~\cite{BFRY} and James,
Lijoi and Pr\"{u}nster~\cite{JLP}~(see also~\cite{JamesGamma}, which
is a~preliminary version of this work) show that the study of mean
functionals is relevant to the analysis of phenomena related to
Bessel and Brownian processes. In fact, the work of James, Lijoi
and Pr\"{u}nster~\cite{JLP} identifies many new explicit examples of
Dirichlet means which have interesting interpretations.

Related to these last points, Lijoi and Regazzini~\cite{Lijoi}
have highlighted a close connection to the theory of generalized
gamma convolutions~(see~\cite{BondBook}). Specifically, it is
known that a rich subclass of random variables having generalized
gamma convolutions~(GGC) distributions may be represented as
%
\begin{equation}
T_{\theta}\stackrel{\mathrm{d}}=G_{\theta}M_{\theta}(F_{X})\stackrel{\mathrm{d}}=\int
_{0}^{\infty}x\gamma_{\theta}(\mathrm{d}x).
\label{gammarep}
\end{equation}
We call these random variables GGC$(\theta,F_{X}).$ In addition, we
see from~(\ref{gammarep}) that $T_{\theta}$ is a random variable
derived from a weighted gamma process and, hence, the calculus
discussed in Lo~\cite{Lo82} and Lo and Weng~\cite{LW} applies. In
general, GGC random variables are an important class of infinitely
divisible random variables whose properties have been extensively
studied by~\cite{BondBook} and others. We note further that
although we have written a GGC$(\theta,F_{X})$ random variable as
$G_{\theta}M_{\theta}(F_{X})$, this representation is not unique
and, in fact, it is quite rare to see $T_{\theta}$ represented in
this way. We will show that one can, in fact, exploit this
non-uniqueness to obtain explicit densities for $T_{\theta}$, even
when it is not so easy to do so for $M_{\theta}(F_{X}).$ While the
representation $G_{\theta}M_{\theta}(F_{X})$ is not unique, it
helps one to understand the relationship between the Laplace
transform of $T_{\theta}$ and the Cauchy--Stieltjes transform of
order $\theta$ of $M_{\theta}(F_{X}),$ which, indeed, characterize
respectively the laws of $T_{\theta}$ and $M_{\theta}(F_{X}).$
Specifically, using the independence property of $G_{\theta}$ and
$M_{\theta}(F_{X})$ leads to, for $\lambda\ge0,$
%
\begin{equation}
\mathbb{E}[\mathrm{e}^{-\lambda T_{\theta}}]=\mathbb{E}\bigl[{\bigl(1+\lambda
M_{\theta}(F_{X})\bigr)}^{-\theta}\bigr]=\mathrm{e}^{-\theta\psi_{F_{X}}(\lambda)}, \label{CS}
\end{equation}
where
%
\begin{equation}
\psi_{F_{X}}(\lambda)=\int_{0}^{\infty}\log(1+\lambda
x)F_{X}(\mathrm{d}x)=\mathbb{E}[\log(1+\lambda X)] \label{Levy}
\end{equation}
is the \textit{L\'{e}vy exponent} of $T_{\theta}.$ We note that
$T_{\theta}$ and $M_{\theta}(F_{X})$ exist if and only if
$\psi_{F_{X}}(\lambda)<\infty$ for $\lambda>0$ (see, e.g.,~\cite{DS}
and~\cite{BondBook}). The expressions
in~(\ref{CS}) equate with the aforementioned identity obtained by
Cifarelli and
Regazzini~\cite{CifarelliRegazzini}.

Despite these interesting results, there is very little work on
the relationship between different mean functionals. Suppose, for
instance, that for each fixed value of $\theta>0,$
$M_{\theta}(F_{X})$ denotes a~Dirichlet mean and $(M_{\theta}(F_{Z_{c}});c>0)$
denotes a collection of Dirichlet mean random variables indexed
by a family of distributions $(F_{Z_{c}};c>0 ).$
One may then ask the following question: for what choices of $X$ and
$Z_{c}$ are
these mean functionals related, and in what sense? In particular,
one may wish to know how their densities are related. The
rationale here is that if such a relationship is established, then the
effort that
one puts forth to obtain results such as the explicit density of
$M_{\theta}(F_{X})$
can be applied to an entire family of Dirichlet means
$(M_{\theta}(F_{Z_{c}});c>0).$ Furthermore, since
Dirichlet means are associated with GGC random variables, this
would establish relationships between a GGC$(\theta,F_{X})$ random
variable and a family of GGC$(\theta,F_{Z_{c}})$ random
variables.
Simple examples are, of course, the choices $Z_{c}=X+c$ and
$Z_{c}=c X,$ which, due to the linearity properties of mean
functionals, result easily in the identities in law
\[
M_{\theta}(F_{X+c})=c+M_{\theta}(F_{X})\quad \mbox{and}\quad
M_{\theta}(F_{c X})=cM_{\theta}(F_{X}).
\]

Naturally, we are going to discuss more complex relationships, but
with the same goal. That is, we will identify non-trivial
relationships so that the often considerable efforts that one
makes in the study of one mean functional $M_{\theta}(F_{X})$ can
then be used to obtain more easily results for other mean
functionals, their corresponding GGC random variables and L\'{e}vy
processes. In this paper, we will describe two such operations
which we elaborate on in the next subsection.

\subsection{Outline and summary of results}\label{sec:outline}
Section~\ref{sec:prelim} reviews some of the existing formulae for
the densities and cdfs of Dirichlet means. In Section~\ref{sec:beta},
we will describe the operation of multiplying a mean functional
$M_{\theta\sigma}(F_{X})$ by an independent beta random variable
with parameters $(\theta\sigma,\theta(1-\sigma)),$ say,
$\beta_{\theta\sigma,\theta(1-\sigma)}$, where $0<\sigma<1.$ We
call this operation \textit{beta scaling}. Theorem~\ref{thm1beta}
shows that the resulting random variable
$\beta_{\theta\sigma,\theta(1-\sigma)}M_{\theta\sigma}(F_{X})$ is
again a mean functional, but now of order $\theta$. In addition,
the GGC$(\theta\sigma,F_{X})$ random variable $G_{\theta
\sigma}M_{\theta\sigma}(F_{X})$ is equivalently a GGC random
variable of order $\theta.$ Now, keeping in mind that tractable
densities of mean functionals of order $\theta=1$ are the easiest
to obtain, Theorem~\ref{thm1beta} shows that by setting
$\theta=1$, the densities of the uncountable collection of random
variables $(\beta_{\sigma,1-\sigma}M_{\sigma}(F_{X}); 0<\sigma\leq
1)$ are all mean functionals of order $\theta=1.$
Theorem~\ref{thm2beta} then shows that efforts used to calculate
the explicit density of any one of these random variables, via the
formulae of~\cite{CifarelliRegazzini}, lead to the explicit
calculation of the densities of all of them. Additionally,
Theorem~\ref{thm2beta} shows that the corresponding GGC random
variables may all be expressed as GGC random variables of order
$\theta=1,$ representable in distribution as
$G_{1}\beta_{\sigma,1-\sigma}M_{\sigma}(F_{X})$. A key point here
is that Theorem~\ref{thm2beta} gives a tractable density for
$\beta_{\sigma,1-\sigma}M_{\sigma}(F_{X})$ without requiring
knowledge of the density of $M_{\sigma}(F_{X}),$ which is usually
expressed in a complicated manner. These results will also yield
some non-obvious integral identities. Furthermore, noting that a
GGC$(\theta,F_{X})$ random variable, $T_{\theta},$ is infinitely
divisible, we associate it with an independent increment process
$(\zeta_{\theta}(t); t\ge0),$ known as a \textit{subordinator}
(a~non-decreasing non-negative L\'{e}vy process), where, for each fixed
$t,$
\[
\mathbb{E}\bigl[\mathrm{e}^{-\lambda
\zeta_{\theta}(t)}\bigr]=\mathbb{E}[\mathrm{e}^{-\lambda
T_{\theta t}}]=\mathrm{e}^{-t\theta\psi_{F_{X}}(\lambda)}.
\]
That is, marginally, $\zeta_{\theta}(1)\stackrel{\mathrm{d}}=T_{\theta}$ and
$\zeta_{\theta}(t)\stackrel{\mathrm{d}}=\zeta_{\theta
t}(1)\stackrel{\mathrm{d}}=T_{\theta t}.$ In addition, for $s<t,$
$\zeta_{\theta}(t)-\zeta_{\theta}(s)\stackrel{\mathrm{d}}=\zeta_{\theta}(t-s)$
is independent of $\zeta_{\theta}(s).$ We say that the process
$(\zeta_{\theta}(t); t\ge0)$ is a \textit{GGC$(\theta,F_{X})$
subordinator}. Proposition~\ref{prop1beta} shows how
Theorems~\ref{thm1beta} and~\ref{thm2beta} can be used to address
the usually difficult problem of explicitly describing the
densities of the finite-dimensional distribution of
a~subordinator~(see \cite{Kingman75}). This has implications in, for
instance, the explicit description of densities of Bayesian
nonparametric prior and posterior models, but is clearly of wider
interest in terms of the distribution theory of infinitely
divisible random variables and associated processes.

In Section~\ref{sec:gamma}, we describe how the operation of
exponentially tilting the density of\break a~GGC$(\theta, F_{X})$ random
variable leads to a relationship between the densities of the mean
functional $M_{\theta}(F_{X})$ and yet another family of mean
functionals. This is summarized in Theorem~\ref{thm1gamma}.
Section~\ref{sec:tiltbeta} then discusses a combination of the two
operations. Proposition~\ref{prop1gamma} describes the density of
beta-scaled and tilted mean functionals of order 1. Using this,
Proposition~\ref{prop2gamma} describes a~method to calculate a key
quantity in the explicit description of the densities and cdfs of
mean functionals. In Section~\ref{sec:example}, we show how the
results in Sections \ref{sec:beta} and \ref{sec:gamma} are used to derive the finite-dimensional
distribution and related quantities for classes of
subordinators suggested by the recent work of James, Lijoi and
Pr\"{u}nster~\cite{JLP} and ~Bertoin \textit{et al.}~\cite{BFRY}.

\subsection{Preliminaries}\label{sec:prelim}
Suppose that $X$ is a positive random variable with
distribution $F_{X}$ and define the function
\[
\Phi_{F_{X}}(t)=\int_{0}^{\infty}\log(|t-x|)\indic(t\neq
x)F_{X}(\mathrm{d}x)=\mathbb{E}[\log(|t-X|)\indic(t\neq X)].
\]
Furthermore, define
\[
\Delta_{\theta}(t|F_{X})=\frac{1}{\curpi}\sin(\curpi\theta
F_{X}(t))\mathrm{e}^{-\theta\Phi_{F_{X}}(t)},
\]
where, using a Lebesgue--Stieltjes integral,
$F_{X}(t)=\int_{0}^{t}F_{X}(\mathrm{d}x).$ Cifarelli and
Regazzini~\cite{CifarelliRegazzini}~(see
also~\cite{CifarelliMelilli}) apply an inversion formula to obtain
the distributional formula for $M_{\theta}(F_{X})$ as follows. For
all $\theta>0$, the cdf can be expressed as
%
\begin{equation}
\label{DPcdf} \int_{0}^{x}{(x-t)}^{\theta-1} \Delta_{\theta
}(t|F_{X})\,\mathrm{d}t ,
\end{equation}
provided that $\theta F_{X}$ possesses no jumps of size
greater than or equal to one. If we let $\xi_{\theta
F_{X}}(\cdot)$ denote the density of $M_{\theta}(F_{X}),$
then it takes its simplest form for $\theta=1$, which is
%
\begin{equation} \label{M1}
\xi_{F_{X}}(x)=\Delta_{1}(x|F_{X})=\frac{1}{\curpi}\sin(\curpi
F_{X}(x))\mathrm{e}^{-\Phi(x)}.
\end{equation}
Density formulae for
$\theta>1$ are described as
%
\begin{equation} \label{M2}
\xi_{\theta F_{X}}(x)= {(\theta-1)}\int_{0}^{x}{(x-t)}^{\theta-2}
\Delta_{\theta}(t|F_{X})\,\mathrm{d}t.
\end{equation}

An expression for the density, which holds for all $\theta>0$, was
recently obtained by James, Lijoi and Pr\"{u}nster~\cite{JLP} as
follows:
%
\begin{equation} \label{generaldensity} \xi_{\theta F_{X}}(x)=
\frac{1}{\curpi}\int_{0}^{x}{(x-t)}^{\theta-1} d_{\theta}(t|F_{X})\,\mathrm{d}t ,
\end{equation}
where
\[
d_{\theta}(t|F_{X})=\frac{\mathrm{d}}{\mathrm{d}t}\sin(\curpi\theta
F_{X}(t))\mathrm{e}^{-\theta\Phi(t)}.
\]
For additional formulae,
see~\cite{CifarelliRegazzini,Regazzini2002,Lijoi}.

\begin{rem} Throughout, for random variables $R$ and $X,$ when we
write the product $RX$, we will assume, unless otherwise
mentioned, that $R$ and $X$ are independent. This convention
will also apply to the multiplication of the special random
variables that are expressed as mean functionals. That is,
the product $M_{\theta}(F_{X})M_{\theta}(F_{Z})$ is
understood to be a product of independent Dirichlet means.
\end{rem}

\begin{rem} Throughout, we will be using the fact that if $R$ is a gamma
random variable, then the independent random variables $R,X,Z$
satisfying $RX\stackrel{\mathrm{d}}=RZ$ imply that $X\stackrel{\mathrm{d}}=Z.$ This is
true because gamma random variables are simplifiable. For the precise
meaning of this term and associated conditions, see Chaumont and
Yor~\cite{Chaumont}, Sections
1.12 and 1.13. This fact also applies to the case where
$R$ is a positive stable random variable.
\end{rem}

\section{Beta scaling}~\label{sec:beta}In this section, we investigate
the simple operation of multiplying a Dirichlet mean functional
$M_{\theta}(F_{X})$ by certain beta random variables. Note, first,
that if $M$ denotes an arbitrary positive random variable with
density $f_{M},$ then, by elementary arguments, it follows that the
random variable $W\stackrel{\mathrm{d}}=\beta_{a,b}M,$ where $\beta_{a,b}$
is beta$(a,b)$ independent of $M,$ has density expressible as
\[
f_{W}(w)=\frac{\Gamma(a+b)}{\Gamma(a)\Gamma(b)}\int
_{0}^{1}f_{M}(w/u)u^{a-2}{(1-u)}^{b-1}\,\mathrm{d}u.
\]
However, it is only in special cases that the density
$f_{W}$ can be expressed in even simpler terms. That is to
say, it is not obvious how to carry out the integration. In
the next results, we show how remarkable simplifications can
be achieved when $M=M_{\theta}(F_{X}),$ in particular, for
the range $0<\theta\leq1,$ and when $\beta_{a,b}$ is a
symmetric beta random variable. First, we will need to
introduce some additional notation. Let $Y_{\sigma}$ denote
a Bernoulli random variable with success probability
$0<\sigma\leq1.$ Then, if $X$ is a random variable with
distribution $F_{X}$, independent of $Y_{\sigma},$ it
follows that the random variable $XY_{\sigma}$ has
distribution
%
\begin{equation}
F_{XY_{\sigma}}(\mathrm{d}x)=\sigma
F_{X}(\mathrm{d}x)+(1-\sigma)\delta_{0}(\mathrm{d}x) \label{pdfp}
\end{equation}
and cdf
%
\begin{equation}
F_{XY_{\sigma}}(x)=\sigma F_{X}(x)+(1-\sigma)\indic(x\ge0)
\label{cdfp}.
\end{equation}

Hence, there exists the mean functional
\[
M_{\theta}(F_{XY_{\sigma}})\stackrel{\mathrm{d}}=\int_{0}^{\infty}y\tilde
{P}_{\theta}(\mathrm{d}y) ,
\]
where $\tilde{P}_{\theta}(\mathrm{d}y)$ denotes a Dirichlet process
with parameters $(\theta, F_{XY_{\sigma}}).$ In addition, we
have, for $x>0,$
%
\begin{equation}
\Phi_{F_{XY_{\sigma}}}(x)=\mathbb{E}[\log(|x-XY_{\sigma}|)\indic
(XY_{\sigma}\neq
x)]=\sigma\Phi_{F_{X}}(x)+(1-\sigma)\log(x) \label{Pid}.
\end{equation}
When $\sigma=1,$ $Y_{\sigma}=1$ and hence
$F_{XY_{1}}(x)=F_{X}(x).$ Let $E_{\sigma}$ denote a set
such that $\mathbb{E}[P_{\theta}(E_{\sigma})]=\sigma.$ Note, now, that
every beta random variable, $\beta_{a,b},$
where $a, b$ are arbitrary positive constants, can be
represented as the simple mean functional
\[
P_{\theta}(E_{\sigma})\stackrel{\mathrm{d}}=\beta_{\theta
\sigma,\theta(1
-\sigma)}\stackrel{\mathrm{d}}=M_{\theta}(F_{Y_{\sigma}}),
\]
by
choosing
\[
\sigma=\frac{a}{a+b}\quad \mbox{and}\quad \theta=a+b.
\]
We note, however, that there are other choices of $F_{X}$
that will also yield beta random variables as mean
functionals. Throughout, we will use the convention that
$\beta_{\theta,0}:=1,$ that is, the case where $\sigma=1.$ We
now present our first result.

\begin{thm}~\label{thm1beta} For $\theta>0$ and $0<\sigma\leq1$, let
$\beta_{\theta\sigma,\theta(1 -\sigma)}$ denote a beta
random variable with parameters $(\theta\sigma,
\theta(1-\sigma))$, independent of the mean functional
$M_{\theta\sigma}(F_{X}).$ Then:
\begin{longlist}[(iii)]
\item[(i)]$\beta_{\theta\sigma,\theta(1-\sigma)}M_{\theta
\sigma}(F_{X})\stackrel{d}=M_{\theta}(F_{XY_{\sigma}});$
\item[(ii)]equivalently, $M_{\theta}(F_{Y_{\sigma}})M_{\theta
\sigma}(F_{X})\stackrel{d}=M_{\theta}(F_{XY_{\sigma}});$
\item[(iii)]$G_{\theta\sigma}M_{\theta\sigma}(F_{X})\stackrel
{d}=G_{\theta}M_{\theta
}(F_{XY_{\sigma}});$
\item[(iv)]that is, GGC$(\theta\sigma,F_{X})=$GGC$(\theta,F_{XY_{\sigma}})$.
\end{longlist}
\end{thm}
\begin{pf}Since $M_{\theta}(F_{Y_{\sigma}})\stackrel{\mathrm{d}}=\beta_{\theta
\sigma,\theta(1-\sigma)}$, statements~(i) and~(ii) are
equivalent. We proceed by first establishing~(iii)
and~(iv). Note that, using~(\ref{Levy}),
\[
\mathbb{E}[\log(1+\lambda XY_{\sigma})]=\sigma
\mathbb{E}[\log(1+\lambda
X)]=\sigma\int_{0}^{\infty}\log(1+\lambda x)F_{X}(\mathrm{d}x).
\]
Hence,
\[
\mathbb{E}\bigl[\mathrm{e}^{-\lambda
G_{\theta}M_{\theta}(F_{XY_{\sigma}})}\bigr]=\mathrm{e}^{-\theta
\sigma\int_{0}^{\infty}\log(1+\lambda
x)F_{X}(\mathrm{d}x)}=\mathbb{E}\bigl[\mathrm{e}^{-\lambda G_{\theta\sigma
}M_{\theta\sigma}(F_{X})}\bigr],
\]
which means that
$G_{\theta}M_{\theta}(F_{XY_{\sigma}})\stackrel{\mathrm{d}}=G_{\theta\sigma}M_{\theta\sigma}(F_{X}),$ establishing statements~(iii)
and~(iv). Now, writing
$G_{\theta\sigma}=G_{\theta}\beta_{\theta
\sigma,\theta(1-\sigma)},$ it follows that
\[
G_{\theta}M_{\theta}(F_{XY_{\sigma}})\stackrel{\mathrm{d}}=G_{\theta}\beta_{\theta
\sigma,\theta(1-\sigma)}M_{\theta\sigma}(F_{X}).
\]
Hence, $
\beta_{\theta
\sigma,\theta(1-\sigma)}M_{\theta\sigma}(F_{X})\stackrel{\mathrm{d}}=M_{\theta
}(F_{XY_{\sigma}}),
$ by the fact that gamma random variables are simplifiable.
\end{pf}
\begin{rem}We note that parts (i) and (ii) of Theorem 2.1 also follow as
consequences of Ethier and Griffiths~\cite{Ethier}, Lemma 1. We
now state their interesting result for clarity.
\begin{lem}[(Ethier and Griffiths~\cite{Ethier})]
Let $\nu_{1}$ and $\nu_{2}$ denote two probability measures. Now,
for $\theta_{1},\theta_{2}>0,$ define the probability measure
\[
\nu_{(\theta_{1},\theta_{2})}(\mathrm{d}x)=\frac{\theta_{1}}{\theta_{1}+\theta
_{2}}\nu_{1}(\mathrm{d}x)+\frac{\theta_{2}}{\theta_{1}+\theta_{2}}\nu_{2}(\mathrm{d}x).
\]
Then, for independent Dirichlet processes $\mu_{1}\sim
\Pi_{\theta_{1},\nu_{1}}$ and $\mu_{2}\sim
\Pi_{\theta_{2},\nu_{2}}$,
\[
\mu_{1,2}(\cdot)\stackrel{d}=\beta_{\theta_{1},\theta_{2}}\mu_{1}(\cdot
)+(1-\beta_{\theta_{1},\theta_{2}})\mu_{2}(\cdot),
\]
where $\mu_{1,2}$ is a Dirichlet process with parameters
$(\theta_{1}+\theta_{2},\nu_{(\theta_{1},\theta_{2})}).$
\end{lem}

Hence, as a general consequence,
\[
M_{\theta_{1}+\theta_{2}}\bigl(\nu_{(\theta_{1},\theta_{2})}\bigr)\stackrel
{d}=\beta_{\theta_{1},\theta_{2}}M_{\theta_{1}}(\nu_{1})
+(1-\beta_{\theta_{1},\theta_{2}})M_{\theta_{2}}(\nu_{2}).
\]
Now, from~(\ref{pdfp}), we see that setting
$\nu_{1}=F_{X},\nu_{2}=\delta_{0}, \theta_{1}=\theta\sigma$ and
$\theta_{2}=\theta(1-\sigma)$ yields statements \textup{(i)} and \textup{(ii)}.
This is because $M_{\theta(1-\sigma)}(\delta_{0})=0.$
\end{rem}

 When $\theta=1$, we obtain results for random variables
$\beta_{\sigma,1-\sigma}M_{\sigma}(F_{X}).$ The symmetric
beta random variables $\beta_{\sigma,1-\sigma}$ arise in a
variety of important contexts and are often referred to as
generalized arcsine laws with density expressible as
\[
\frac{\sin(\curpi
\sigma)}{\curpi}u^{\sigma-1}{(1-u)}^{-\sigma}\qquad{\mbox{for
}}0<u<1.
\]
Now, using (\ref{pdfp}) and (\ref{cdfp}), let $
\mathcal{C}(F_{X})=\{x:F_{X}(x)>0\}$. Then, for $x>0,$
%
\begin{equation}
\label{sinp} \sin(\curpi
F_{XY_{\sigma}}(x))=\cases{
\sin\bigl(\curpi\sigma[1-F_{X}(x)]\bigr), &\quad  $\mbox{if }
x\in\mathcal{C}(F_{X}),$\cr
\sin(\curpi(1-\sigma)),&\quad  $\mbox{if } x\notin\mathcal{C}(F_{X}).$}
\end{equation}
Also, note that $\sin(\curpi[1-F_{X}(x)])=\sin(\curpi F_{X}(x)).$
The next result yields another surprising property of these
random variables.

\begin{thm}~\label{thm2beta}Consider the setting in Theorem~\ref
{thm1beta}. Then, when $\theta=1$, it follows that for each fixed
$0<\sigma\leq1,$ the random variable
$M_{1}(F_{XY_{\sigma}})\stackrel{d}=\beta_{\sigma,1-\sigma}M_{\sigma}(F_{X})$
has density
%
\begin{equation}
\xi_{F_{XY_{\sigma}}}(x)=\frac{x^{\sigma-1}}{\curpi}\sin(\curpi
F_{XY_{\sigma}}(x))\mathrm{e}^{-\sigma\Phi_{F_{X}}(x)}\qquad\mbox{for } x>0, \label{betaid}
\end{equation}
specified by~(\ref{sinp}). Since
GGC$(\sigma,F_{X})=$GGC$(1,F_{XY_{\sigma}}),$ this implies
that the random variable
$G_{\sigma}M_{\sigma}(F_{X})\stackrel{d}=G_{1}M_{1}(F_{XY_{\sigma}})$
has density
%
\begin{equation}
g_{\sigma,F_{X}}(x)=\frac{1}{\curpi}\int_{0}^{\infty}\mathrm{e}^{-{x}/{y}}y^{\sigma-2}\sin(\curpi
F_{XY_{\sigma}}(y))\mathrm{e}^{-\sigma\Phi_{F_{X}}(y)}\,\mathrm{d}y .
\label{mix}
\end{equation}
\end{thm}
\begin{pf}Since
$M_{1}(F_{XY_{\sigma}})\stackrel{\mathrm{d}}=\beta_{\sigma,1-\sigma}M_{\sigma}(F_{X}),$
the density is of the form~(\ref{M1}) for each fixed
$\sigma\in(0,1].$ Furthermore, we use the identity
in~(\ref{Pid}).
\end{pf}
\begin{rem}
It is worthwhile to mention that transforming to the random
variable
$1/\beta_{\sigma,1-\sigma}$,~(\ref{betaid}) is equivalent
to the otherwise non-obvious integral identity
\[
\frac{\sin(\curpi
\sigma)}{\curpi}\int_{1}^{\infty}\frac{\xi_{\sigma
F_{X}}(xy)}{{(y-1)}^{\sigma}}\,\mathrm{d}y
=\frac{x^{\sigma-1}}{\curpi}\sin(\curpi F_{XY_{\sigma}}(x))\mathrm{e}^{-\sigma\Phi(x)}.
\]
This leads to interesting results when the density
$\xi_{\sigma F_{X}}(x)$ has a known form. On the other
hand, we see that one does not need the explicit density of
$M_{\sigma}(F_{X})$ to obtain the density of
$M_{1}(F_{XY_{\sigma}})\stackrel{\mathrm{d}}=\beta_{\sigma,1-\sigma}M_{\sigma}(F_{X}).$
In fact, owing to our goal of yielding simple densities for
many Dirichlet means from one mean, we see that the effort
to calculate the density of $M_{1}(F_{XY_{\sigma}})$ for
each $0< \sigma\leq1$ is no more than what is needed to
calculate the density of $M_{1}(F_{X}).$
\end{rem}

We now see how this translates into the usually difficult problem
of explicitly describing the density of the finite-dimensional
distribution of a subordinator. In the next result, we write, for
some set $C,$
\[
\zeta_{\theta}(C):=\int_{0}^{\infty}\indic(s\in
C)\zeta_{\theta}(\mathrm{d}s).
\]

\begin{prop}~\label{prop1beta}Let $(\zeta_{\theta}(t);t\leq1/\theta)$
denote a GGC$(\theta,F_{X})$ subordinator on $[0,1/\theta].$ Furthermore,
let $(C_{1}, \ldots, C_{k})$ denote an arbitrary disjoint
partition of the interval $[0,1/\theta].$ The
finite-dimensional distribution
$(\zeta_{\theta}(C_{1}),\ldots,\zeta_{\theta}(C_{k}))$ then has
a joint density
%
\begin{equation} \prod_{i=1}^{k}g_{\sigma_{i},F_{X}}(x_{i}),
\label{fidi1}
\end{equation}
where each $\sigma_{i}=\theta|C_{i}|>0$ and
$\sum_{i=1}^{k}\sigma_{i}=1.$ The density
$g_{\sigma_{i},F_{X}}$ is given by~(\ref{mix}). That is,
$\zeta_{\theta}(C_{i})\stackrel{d}=G_{1}M_{1}(F_{XY_{\sigma_{i}}})$
and these are independent for $i=1,\ldots, k,$ where
$M_{1}(F_{XY_{\sigma_{i}}})\stackrel{d}=\beta_{\sigma_{i},1-\sigma
_{i}}M_{\sigma_{i}}(F_{X})$
has density
\[
\frac{1}{\curpi}x^{\sigma_{i}-1}\sin(\curpi
F_{XY_{\sigma_{i}}}(x))\mathrm{e}^{-\sigma_{i}\Phi_{F_{X}}(x)}.
\]
\end{prop}

\begin{pf} First, since $(C_{1},\ldots,C_{k})$
partitions the interval $[0,1/\theta],$ it follows that
their sizes satisfy $0<|C_{k}|\leq1/\theta$ and
$\sum_{i=1}^{k}|C_{k}|=1/\theta.$ Since $\zeta_{\theta}$ is
a subordinator, the independence of the
$\zeta_{\theta}(C_{i})$ is a consequence of its independent
increment property. In fact, these are essentially
equivalent statements. Hence, we can isolate each
$\zeta_{\theta}(C_{i}).$ It follows that for each $i$, the
Laplace transform is given by
\[
\mathbb{E}\bigl[\mathrm{e}^{-\lambda
\zeta_{\theta}(C_{i})}\bigr]=\mathrm{e}^{-\theta
|C_{i}|\psi_{F_{X}}(\lambda)}=\mathrm{e}^{-\sigma_{i}\psi_{F_{X}}(\lambda)},
\]
which shows that each $\zeta_{\theta}(C_{i})$ is
GGC$(\sigma_{i},F_{X})$ for $0<\sigma_{i}\leq1.$ Hence, the
result follows from Theorem~\ref{thm2beta}.
\end{pf}

\section{Exponential tilting/Esscher transform}\label{sec:gamma}
In this section, we describe how the operation of \textit{exponential
tilting} of the density of a~GGC$(\theta,F_{X})$ random variable
leads to a non-trivial relationship between a mean functional
determined by $F_{X}$ and $\theta,$ and an entire family of mean
functionals indexed by an arbitrary constant $c>0.$ Additionally,
this will identify a non-obvious relationship between two classes
of mean functionals. Exponential tilting is merely a convenient phrase
for the operation of applying an exponential change of measure to
a density or more general measure. In mathematical finance and
other applications, it is known as an \textit{Esscher transform}
and is a key tool for option pricing. We mention that there is
much known about exponential tilting of infinitely divisible
random variables and, in fact, Bondesson~\cite{BondBook}, Example
3.2.5, explicitly discusses the case of GGC random
variables, albeit not in the way we shall describe it. In
addition, examining the gamma representation in~(\ref{gammarep}),
one can see a relationship to Lo and Weng~\cite{LW}, Proposition
3.1 (see also K\"{u}chler and Sorensen~\cite{Kuc} and
James~\cite{James05}, Proposition 2.1), for results on exponential
tilting of L\'{e}vy processes). However, our focus here is on the
properties of related mean functionals, which leads to genuinely
new insights.

Before we elaborate on this, we describe generically what we mean
by exponential tilting. Suppose that $T$ denotes an arbitrary
positive random variable with density, say, $f_{T}.$ It follows
that for each positive $c,$ the random variable $cT$ is
well defined and has density
\[
\frac{1}{c}f_{T}(t/c).
\]
Exponential tilting refers to the exponential change of measure
resulting in a random variable, say $\tilde{T}_{c},$ defined by
the density
\[
f_{\tilde{T}_{c}}(t)=\frac{\mathrm{e}^{-t}({1}/{c})f_{T}(t/c)}{\mathbb{E}[\mathrm{e}^{-cT}]}.
\]
Thus, from the random variable $T$, one gets a family of random
variables $(\tilde{T}_{c};c>0).$ Obviously, the density for each
$\tilde{T_{c}}$ does not differ much. However, something
interesting happens when $T$ is a scale mixture of a gamma random
variables, that is, $T=G_{\theta}M$ for some random positive random
variable $M$ independent of $G_{\theta}.$ In that case, one can
show, see \cite{JamesGamma}, that $T_{c}=G_{\theta}\tilde{M}_{c}$,
where $\tilde{M}_{c}$ is sufficiently distinct for each value of
$c.$ We demonstrate this for the case where $M=M_{\theta}(F_{X}).$

First, note that, obviously,
$cM_{\theta}(F_{X})=M_{\theta}(F_{cX})$ for each $c>0,$ which, in
itself, is not a~very interesting transformation. Now, setting
$T_{\theta}=G_{\theta}M_{\theta}(F_{X})$ with density denoted
$g_{\theta,F_{X}}$, the corresponding
random variable $\tilde{T}_{\theta,c}$
resulting from exponential tilting has density
%
\begin{equation}
\mathrm{e}^{-t}(1/c)g_{\theta,F_{X}}(t/c)\mathrm{e}^{\theta\psi_{F_{X}}(c)}
\label{tilt1}
\end{equation}
and Laplace transform
%
\begin{equation}
\frac{\mathbb{E}[\mathrm{e}^{-c(1+\lambda)G_{\theta}M_{\theta}(F_{X})}]}{\mathbb{E}[\mathrm{e}^{-cG_{\theta}M_{\theta}(F_{X})}]}=\mathrm{e}^{-\theta[\psi_{F_{X}}(c(1+\lambda))-\psi_{F_{X}}(c)]}.
\label{tilt2}
\end{equation}
Now, for each $c>0,$ define the random variable
\[
A_{c}\stackrel{\mathrm{d}}=\frac{cX}{(cX+1)}.
\]
That is, the cdf of the random
variable $A_{c}$ can be expressed as
\[
F_{A_{c}}(y)=F_{X}\biggl(\frac{y}{c(1-y)}\biggr)\qquad{\mbox{for
}}0<y<1.
\]
In the next theorem, we will show that $M_{\theta}(F_{X})$ relates
to the family of mean functionals $(M_{\theta}(F_{A_{c}});c>0)$
by the tilting operation described above. Moreover, we will
describe the relationship between their densities.

\begin{thm}\label{thm1gamma} Suppose that $X$ has distribution $F_{X}$
and for each $c>0,$ $A_{c}\stackrel{d}=cX/(cX+1)$ is a random variable
with distribution $F_{A_{c}}.$ For each $\theta>0$, let
$T_{\theta}=G_{\theta}M_{\theta}(F_{X})$ denote a GGC
$(\theta,F_{X})$ random variable having density
$g_{\theta,F_{X}}.$ Let $\tilde{T}_{\theta,c}$ denote a random
variable with density and Laplace transform described by
(\ref{tilt1}) and~(\ref{tilt2}), respectively.
$\tilde{T}_{\theta,c}$ is then a GGC$(\theta,F_{A_{c}})$ random
variable and hence representable as
$G_{\theta}M_{\theta}(F_{A_{c}}).$ Furthermore, the following
relationships exist between the densities of the mean functionals
$M_{\theta}(F_{X})$ and $M_{\theta}(F_{A_{c}})$:
\begin{enumerate}[(ii)]
\item[(i)]supposing that the density of $M_{\theta}(F_{X})$, say $\xi
_{\theta F_{X}}$,
is known, then the density of $M_{\theta}(F_{A_{c}})$ is
expressible as
\[
\xi_{\theta F_{A_{c}}}(y)=\frac{1}{c}\mathrm{e}^{\theta\psi_{F_{X}}(c)}{{(1-y)}^{\theta-2}}\xi_{\theta
F_{X}}\biggl(\frac{y}{c(1-y)}\biggr)
\]
for $0<y<1;$
\item[(ii)]conversely, if the density of $M_{\theta}(F_{A_{c}})$,
$\xi_{\theta F_{A_{c}}}(y),$ is known, then the density of
$M_{\theta}(F_{X})$ is given by
\[
\xi_{\theta F_{X}}(x)={(1+x)}^{\theta-2}\xi_{\theta
F_{A_{1}}}\biggl(\frac{x}{1+x}\biggr) \mathrm{e}^{-\theta\psi_{F_{X}}(1)}.
\]
\end{enumerate}
\end{thm}
\begin{pf}We proceed by first examining the L\'{e}vy exponent
$[\psi_{F_{X}}(c(1+\lambda))-\psi_{F_{X}}(c)]$ associated with
$\tilde{T}_{\theta,c}$ as described in~(\ref{tilt2}). Note that
\[
\psi_{F_{X}}\bigl(c(1+\lambda)\bigr)=\int_{0}^{\infty}\log\bigl(1+c(1+\lambda)x\bigr)F_{X}(\mathrm{d}x)
\]
and $\psi_{F_{X}}(c)$ is of the same form with $\lambda=0.$ Hence,
isolating the logarithmic terms, we can focus on the difference
\[
\log\bigl(1+c(1+\lambda)x\bigr)-\log(1+cx).
\]
This is equivalent to
\[
\log\biggl(1+\frac{cx}{1+cx}\lambda\biggr)=
\log\biggl(\frac{1}{1+cx}+\frac{cx}{1+cx}(1+\lambda)\biggr),
\]
showing that $\tilde{T}_{\theta,c}$ is GGC$(\theta,F_{A_{c}}).$
This fact can also be deduced from Proposition 3.1 in Lo and
Weng~\cite{LW}. The next step is to identify the density of
$M_{\theta}(F_{A_c})$ in terms of the density of
$M_{\theta}(F_{X}).$ Using the fact that
$T_{\theta}=G_{\theta}M_{\theta}(F_{X}),$ one may write the
density of $T_{\theta}$ in terms of a~gamma mixture as
\[
g_{\theta,F_{X}}(t)=\frac{t^{\theta-1}}{\Gamma(\theta)}\int_{0}^{\infty
}\mathrm{e}^{-t/m}m^{-\theta}\xi_{\theta F_{X}}(m)\,\mathrm{d}m.
\]
Hence, rearranging terms in~(\ref{tilt1}), it follows that the
density of $\tilde{T}_{\theta,c}$ can be written as
\[
\mathrm{e}^{\theta
\psi_{F_{X}}(c)}\frac{t^{\theta-1}}{\Gamma(\theta)}\int_{0}^{\infty
}\mathrm{e}^{-t{(cm+1)}/{(cm)}}{(cm)}^{-\theta}\xi_{\theta F_{X}}(m)\,\mathrm{d}m.
\]
Now, further algebraic manipulation makes this look like a mixture
of a gamma$(\theta,1)$ random variables, as follows,
\[
\frac{t^{\theta-1}}{\Gamma(\theta)}\int_{0}^{\infty}\mathrm{e}^{-t{(cm+1)}/{(cm)}}{\biggl[\frac{cm+1}{cm}\biggr]}^{\theta}\frac
{\mathrm{e}^{\theta\psi_{F_{X}}(c)}\xi_{\theta
F_{X}}(m)}{{(1+cm)}^{\theta}}\,\mathrm{d}m.
\]
Hence, it is evident that $M_{\theta}(F_{A_{c}})$ has the same
distribution as a random variable $cM/(cM+1)$, where $M$ has
density
\[
{\mathrm{e}^{\theta\psi_{F_{X}}(c)}{(1+cm)}^{-\theta}\xi_{\theta
F_{X}}(m)}.
\]
Thus, statements~(i) and~(ii) follow.
\end{pf}

\subsection{Tilting and beta scaling}\label{sec:tiltbeta} This section
describes what happens when
one applies the exponential tilting operation
relative to a mean functional resulting from beta scaling. Recall
that the tilting operation applied to
$G_{\theta}M_{\theta}(F_{X})$ described in the previous section
sets up a relationship between $M_{\theta}(F_{X})$ and
$M_{\theta}(F_{Ac}).$ Consider the random variable $\beta_{\theta
\sigma,\theta(1-\sigma)}M_{\theta
\sigma}(F_{X})\stackrel{\mathrm{d}}=M_{\theta}(F_{XY_{\sigma}}).$ Then,
tilting $G_{\theta}M_{\theta}(F_{XY_{\sigma}})$ as in the previous
section leads to the random variable
$G_{\theta}M_{\theta}(F_{cXY_{\sigma}/(cXY_{\sigma}+1)})$ and
hence relates
\[
\beta_{\theta\sigma,\theta(1-\sigma)}M_{\theta
\sigma}(F_{X})\stackrel{\mathrm{d}}=M_{\theta}(F_{XY_{\sigma}})
\]
to the
Dirichlet mean of order $\theta,$
\[
M_{\theta}\bigl(F_{cXY_{\sigma}/(cXY_{\sigma}+1)}\bigr).
\]
Now, letting $F_{A_{c}Y_{\sigma}}$ denote the distribution
of $A_{c}Y_{\sigma},$ one has
\[
A_{c}Y_{\sigma}\stackrel{\mathrm{d}}=\frac{cXY_{\sigma}}{(cXY_{\sigma}+1)}
\]
and hence
%
\begin{equation}
M_{\theta}\bigl(F_{cXY_{\sigma}/(cXY_{\sigma}+1)}\bigr)\stackrel{\mathrm{d}}=M_{\theta}
(F_{A_{c}Y_{\sigma}})\stackrel{\mathrm{d}}=\beta_{\theta\sigma,\theta(1-\sigma
)}M_{\theta\sigma}(F_{A_{c}}).
\end{equation}
In a way, this shows that the order of beta scaling and
tilting can be interchanged. We now derive a result for the
cases of
$M_{1}(F_{XY_{\sigma}})=\beta_{\sigma,1-\sigma}M_{\sigma}(F_{X})$
and
$M_{1}(F_{A_{c}Y_{\sigma}})=\beta_{\sigma,1-\sigma}M_{\sigma}(F_{A_{c}}),$
related by the tilting operation described above. Combining
Theorem~\ref{thm2beta} with Theorem~\ref{thm1gamma} leads
to the following result.

\begin{prop}\label{prop1gamma} For each $0<\sigma\leq1,$ the
random variables
$M_{1}(F_{XY_{\sigma}})=\beta_{\sigma,1-\sigma}M_{\sigma}(F_{X})$
and
$M_{1}(F_{A_{c}Y_{\sigma}})=\beta_{\sigma,1-\sigma}M_{\sigma}(F_{A_{c}})$
satisfy the following:
\begin{enumerate}[(ii)]
\item[(i)]the density of $M_{1}(F_{A_{c}Y_{\sigma}})$ is expressible as
\[
\xi_{F_{A_{c}Y_{\sigma}}}(y)=\frac{\mathrm{e}^{\sigma\psi_{F_{X}}(c)}y^{\sigma-1}}{\curpi
c^{\sigma}{(1-y)}^{\sigma}}\sin\biggl[\curpi
F_{XY_{\sigma}}\biggl(\frac{y}{c(1-y)}\biggr)\biggr]\mathrm{e}^{-\sigma\Phi_{F_{X}}({y}/{(c(1-y))})}
\]
for $0<y<1$;
\item[(ii)]conversely, the density of
$M_{1}(F_{XY_{\sigma}})$ is given by
\[
\xi_{F_{XY_{\sigma}}}(x)=\frac{\mathrm{e}^{-\sigma\psi_{F_{X}}(1)}
x^{\sigma-1}}{\curpi{(1+x)}}\sin\biggl[\curpi
F_{A_{1}Y_{\sigma}}\biggl(\frac{x}{1+x}\biggr)\biggr]\mathrm{e}^{-\sigma\Phi_{F_{A_{1}}}({x}/{(1+x)})}.
\]
\end{enumerate}
\end{prop}
\begin{pf}For clarity, statement (i) is obtained by first using
Theorem~\ref{thm1gamma}, which gives
\[
\xi_{F_{A_{c}Y_{\sigma}}}(y)=\frac{1}{c}\mathrm{e}^{\psi_{F_{XY_{\sigma}}}(c)}{{(1-y)}^{-1}}\xi_{ F_{XY_
{\sigma}}}\biggl(\frac{y}{c(1-y)}\biggr)
\]
for $0<y<1.$ It then
remains to substitute the form of the density~(\ref{betaid}) given
in Theorem~\ref{thm2beta}. Statement (ii) proceeds in the same
way, using~(\ref{mix}).
\end{pf}

Note that even if one can calculate $\Phi_{F_{A_{c}}}$ for
some fixed value of $c,$ it may not be so obvious how to
calculate it for another value. The previous results allow
us to relate their calculation to that of $\Phi_{F_{X}}$, as
described next.

\begin{prop}\label{prop2gamma} Set $A_{c}=cX/(cX+1)$ and define $\Phi
_{F_{A_{c}}}(y)=\mathbb{E}[\log(|y-A_{c}|)\indic(A_{c}\neq
y)].$ Then, for $0<y<1,$
\[
\Phi_{F_{A_{c}}}(y)=\Phi_{F_{X}}\biggl(\frac{y}{c(1-y)}\biggr)
-\psi_{F_{X}}(c)+\log\bigl(c(1-y)\bigr).
\]
\end{prop}

\begin{pf}The result can be deduced by using Proposition~\ref
{prop1gamma} in
the case $\sigma=1.$ First, note that $\sin(\curpi
F_{X}(\frac{y}{c(1-y)}))=\sin(\curpi F_{A_{c}}(y)).$ Now,
equating the form of the density of $M_{1}(F_{A_{c}})$
given by~(\ref{M1}) with the expression given in
Proposition~\ref{prop1gamma}, it follows that
\[
\mathrm{e}^{-\Phi_{F_{A_{c}}}(y)}=\frac{\mathrm{e}^{\psi_{F_{X}}(c)}}{ c{(1-y)}}\mathrm{e}^{-\Phi_{F_{X}}({y}/{(c(1-y))})},
\]
which yields the result.
\end{pf}

\begin{rem} We point out that if $G_{\kappa}$ represents a gamma
random variable for $\kappa\neq\theta$, independent of
$M_{\theta}(F_{X}),$ then it is not necessarily true that
$G_{\kappa}M_{\theta}(F_{X})$ is a GGC random variable. For this
to be true, $M_{\theta}(F_{X})$ would need to be equivalent in
distribution to some $M_{\kappa}(F_{R}).$ In that case, our
results above would be applied for a GGC$(\kappa, F_{R})$ model.
\end{rem}

\section{Examples}\label{sec:example}
In this section, we will demonstrate how our results in
Sections~\ref{sec:beta} and~\ref{sec:gamma} can be applied to
extend results for two random processes recently studied in the
literature. The first involves a class of GGC subordinators that
can be derived from a random mean of a two-parameter Poisson--Dirichlet
process with a uniform base measure, which was studied
as a special case in James, Lijoi and Pr\"{u}nster~\cite{JLP}; see Pitman and
Yor~\cite{PY97} for more details of the two parameter Poisson--Dirichlet distribution. The
second involves a class of processes recently studied in Bertoin
\textit{et al.}~\cite{BFRY}; see also
Maejima~\cite{Maejima} for some discussion of this process. A
key component will be the ability to obtain an explicit expression
for the respective $\Phi_{F_{X}}.$ In the first example, we do not
have much explicit information on the relevant density,
$\xi_{\theta F_{X}};$ however, we can rely on a general theorem of
James, Lijoi and Pr\"{u}nster~\cite{JLP} to obtain $\Phi_{F_{X}}.$
In the second case of the models discussed in Bertoin \textit{et
al.}~\cite{BFRY}, this theorem apparently does not
apply. However, we will be able to use an explicit form of the
density, obtained for a particular value of $\theta$ by Bertoin \textit
{et al.}~\cite{BFRY}, to obtain $\Phi_{F_{X}}.$

As we shall show, both of these processes are connected to a random
variable $Z_{\alpha}$, whose properties we now describe. For
$0<\alpha<1,$ let $S_{\alpha}$ denote a positive $\alpha$-stable
random variable specified by its Laplace transform
\[
\mathbb{E}[\mathrm{e}^{-\lambda S_{\alpha}}]=\mathrm{e}^{-\lambda^{\alpha}}.
\]
In addition, define
\[
Z_{\alpha}={\biggl(\frac{S_{\alpha}}{S'_{\alpha}}\biggr)}^{\alpha},
\]
where $S'_{\alpha}$ is independent of $S_{\alpha}$ and has the
same distribution. The density of this random variable was
obtained by Lamperti~\cite{Lamperti}~(see also Chaumont and
Yor~\cite{Chaumont}, Exercise 4.2.1) and has the remarkably simple
form
\[
f_{Z_{\alpha}}(y)=\frac{\sin(\curpi\alpha)}{\curpi
\alpha}\frac{1}{y^{2}+2y\cos(\curpi\alpha)+1}\qquad{\mbox{for }}y>0.
\]
Furthermore (see Fujita and Yor~\cite{FY} and
(James~\cite{JamesLinnik}, Proposition 2.1), it follows that the cdf
of~$Z_{\alpha}$ satisfies, for $z>0,$

\begin{eqnarray*}
F_{Z_{\alpha}}(1/{z})&=&1-\frac{1}{\curpi\alpha}\cot^{-1}\biggl(\frac{\cos
(\curpi
\alpha)+1/z}{\sin(\curpi\alpha)}\biggr)\\
&=&\frac{1}{\curpi\alpha}\cot^{-1}\biggl(\frac{\cos(\curpi
\alpha)+z}{\sin(\curpi\alpha)}\biggr)\\
&=& 1-F_{Z_{\alpha}}({z}),
\end{eqnarray*}
%
\begin{equation}
\sin(\curpi\alpha F_{Z_{\alpha}}(z))=z\sin\bigl(\curpi
\alpha\bigl(1-F_{Z_{\alpha}}(z)\bigr)\bigr)=\frac{z\sin(\curpi\alpha)}{{[z^{2}+2z\cos(\curpi
\alpha)+1]}^{1/2}} \label{sinid}
\end{equation}
and
%
\begin{eqnarray}
\sin\bigl(2\curpi\alpha[1-F_{Z_{\alpha}}(z)]\bigr)&=&\frac{\sin(2\curpi\alpha)+2z\sin
(\curpi
\alpha)}{1+2z\cos(\curpi\alpha)+z^{2}}
\nonumber
\\[-8pt]
\\[-8pt]
\nonumber
&=&\frac{2\sin(\curpi
\alpha)[\cos(\curpi\alpha)+z]}{1+2z\cos(\curpi\alpha)+z^{2}} .
\label{idd}
\end{eqnarray}
When $\alpha=1/2,$
\[
\sin\bigl(\curpi[1-F_{Z_{1/2}}(z)]\bigr)=\frac{z}{z^{2}+1}.
\]
%
\subsection{Subordinators derived from an example in James, Lijoi and
Pr\"{u}nster} For $0<\alpha<1$ and
$\theta>-\alpha,$ we define the special case of a two-parameter
Poisson--Dirichlet random probability measures as
\[
\tilde{P}_{\alpha,\theta}(\cdot)=\sum_{k=1}^{\infty}V_{k}\prod
_{i=1}^{k-1}(1-V_{i})\delta_{U_{k}}(\cdot),
\]
where $U_{k}$ are i.i.d.~Uniform[0,1] random variables and the
$V_{k}$ are a sequence of independent
$\beta_{\alpha,\theta+k\alpha}$ random variables, independent of
$(U_{k})$. So, in particular, these random variables satisfy
$\mathbb{E}[\tilde{P}_{\alpha,\theta}(\cdot)]=F_{U}(\cdot)$, where
$U$ denotes a Uniform$[0,1]$ random variable. In addition,
$\tilde{P}_{0,\theta}$ is a~Dirichlet process. Then, consider the
random means given as
\[
\tilde{M}_{\alpha,\theta}(F_{U}):=\mathbb{U}_{\alpha,\theta}=\sum
_{k=1}^{\infty}U_{k}V_{k}\prod_{i=1}^{k-1}(1-V_{i})
=\int_{0}^{1}u\tilde{P}_{\alpha,\theta}(\mathrm{d}u).
\]
The $\mathbb{U}_{\alpha,\theta}$ represent a special
case of random variables representable as mean functionals of the
class of two-parameter $(\alpha,\theta)$ Poisson--Dirichlet random
probability measures -- that is to say, random variables
$\tilde{M}_{\alpha,\theta}(F_{X})$, where $F_{X}$ is a general
distribution. An extensive study of this larger class was
conducted by James, Lijoi and Pr\"{u}nster~\cite{JLP}. In regards
to $\mathbb{U}_{\alpha,\theta}$, they show that
$\mathbb{U}_{\alpha,0}$ has an explicit density
\[
\frac{\sin(\curpi\alpha)}{\alpha\curpi}\frac{(\alpha+1)t^{\alpha
}{(1-t)}^{\alpha}}{{[t^{2\alpha+2}+2t^{\alpha+1}
{(1-t)}^{\alpha+1}\cos(\curpi\alpha)+{(1-t)}^{2\alpha+2}]}}.
\]
Furthermore, from James, Lijoi and Pr\"{u}nster~\cite{JLP}, Theorem
2.1, for $\theta>0,$
\[
\mathbb{U}_{\alpha,\theta}\stackrel{\mathrm{d}}=M_{\theta}(F_{\mathbb{U}_{\alpha,0}}).
\]
This implies that
\[
G_{\theta}\mathbb{U}_{\alpha,\theta}\stackrel{\mathrm{d}}=G_{\theta}M_{\theta
}(F_{\mathbb{U}_{\alpha,0}})
\]
are GGC$(\theta,F_{\mathbb{U}_{\alpha,0}}).$ Now, from Vershik, Yor
and Tslevich~\cite{Vershik} (see also James, Lijoi and
Pr\"{u}nster~\cite{JLP}, equation (16)), it follows that
\begin{eqnarray*}
\mathbb{E}[\mathrm{e}^{-\lambda
G_{\theta}\mathbb{U}_{\alpha,\theta}}]&=&{\biggl(\frac{\lambda
(\alpha+1)}{{(\lambda+1)}^{\alpha+1}-1}\biggr)}^{{\theta}/{\alpha}}
\\&=&\exp\bigl(-\theta\mathbb{E}[\log(1+\lambda\mathbb{U}_{\alpha,0})]\bigr),
\end{eqnarray*}
where this expression follows from the generalized Stieltjes
transform of order $-\alpha$ of a Uniform[0,1] random variable,
\[
\mathbb{E}[(1+\lambda U)^{\alpha}]=\int_{0}^{1}{(1+\lambda
x)}^{\alpha}\,\mathrm{d}x=\frac{{(\lambda+1)}^{\alpha+1}-1}{\lambda
(\alpha+1)}.
\]

A description of the densities of $\mathbb{U}_{\alpha,\theta}$ for
$\theta>-\alpha$ is available from the results of~\cite{JLP}.
However, with the exceptions of $\mathbb{U}_{\alpha,1}$ and
$\mathbb{U}_{\alpha,1-\alpha},$ their densities are generally
expressed in terms of integrals with respect to functions that
possibly take on negative values. Here, by focusing instead on
random variables
$\beta_{\theta,1-\theta}\mathbb{U}_{\alpha,\theta}$ for
$0<\theta<1,$ we can utilize the results in James, Lijoi and
Pr\"{u}nster~\cite{JLP} to obtain explicit expressions for their
densities and the corresponding
GGC$(\theta,F_{\mathbb{U}_{\alpha,0}})$ random variables.

In particular, we will obtain explicit descriptions for the
finite-dimensional distribution of a~GGC$(\alpha,F_{\mathbb{U}_{\alpha,0}})$, say
$(\Upsilon_{\alpha}(t),t>0)$, subordinator, where
$\Upsilon_{\alpha}(1)\stackrel{\mathrm{d}}=G_{\alpha}\mathbb{U}_{\alpha,\alpha}$
and hence
\[
\mathbb{E}\bigl[\mathrm{e}^{-\lambda
\Upsilon_{\alpha}(1)}\bigr]=\frac{\lambda
(\alpha+1)}{{(\lambda+1)}^{\alpha+1}-1}.
\]

Although not immediately obvious, one can show that
\[
\mathbb{U}_{\alpha,0}\stackrel{\mathrm{d}}=\frac{Z^{1/(\alpha+1)}_{\alpha
}}{Z^{1/(\alpha+1)}_{\alpha}+1}
{\quad\mbox{and}\quad \mbox{hence }}F_{\mathbb{U}_{\alpha,0}}(t)=F_{Z_{\alpha}}\biggl({\biggl(\frac{t}{1-t}\biggr)}^{\alpha+1}\biggr).
\]
From this, due to the tilting relationship discussed in Section 3,
we see that we can also obtain results for the
GGC$(\alpha,F_{Z^{1/(\alpha+1)}_{\alpha}})$ subordinator, say
$(\Upsilon^{\ddagger}_{\alpha}(t),t>0).$ To the best of our knowledge, this
process and its mean functionals
$M_{\theta}(F_{Z^{1/(\alpha+1}_{\alpha}})$ have not been studied.
Now, from James, Lijoi and Pr\"{u}nster~\cite{JLP}, Theorem 5.2(iii),
it follows that
%
\begin{equation}
\mathrm{e}^{-\Phi_{F_{\mathbb{U}_{\alpha,0}}}(t)}=\frac{{(\alpha+1)}^{1/\alpha
}}{{[t^{2\alpha+2}+2t^{\alpha+1}
{(1-t)}^{\alpha+1}\cos(\curpi
\alpha)+{(1-t)}^{2\alpha+2}]}^{{1}/{(2\alpha)}}} .
\label{JLP}
\end{equation}

This, combined with our results, leads to an explicit description
of the finite-dimensional distribution of the relevant
subordinators.
\begin{thm}Consider the GGC$(\alpha,F_{\mathbb{U}_{\alpha,0}})$ subordinator
$(\Upsilon_{\alpha}(t),t\leq1/\alpha)$ and the\break
GGC$(\alpha,F_{Z^{1/(\alpha+1)}_{\alpha}})$ subordinator
$(\Upsilon^{\ddagger}_{\alpha}(t),t\leq1/\alpha).$ Let
$(C_{1},\ldots,C_{k})$ denote an arbitrary disjoint partition of
the interval $(0,1/\alpha]$ with lengths $|C_{i}|$ and set
$\sigma_{i}=\alpha|C_{i}|$ for $i=1,\ldots,k.$ The following
results then hold:
\begin{enumerate}[(ii)]
\item[(i)] The finite dimensional
distribution of
$(\Upsilon_{\alpha}(C_{1}),\ldots,\Upsilon_{\alpha}(C_{k}))$ is
such that each $\Upsilon_{\alpha}(C_{i})$ is independent and has
distribution
\[
\Upsilon_{\alpha}(C_{i})\stackrel{d}=G_{1}M_{1}(F_{Y_{\sigma_{i}}\mathbb
{U}_{\alpha,0}}),
\]
where
$M_{1}(F_{Y_{\sigma_{i}}\mathbb{U}_{\alpha,0}})\stackrel{d}=\beta_{\sigma
_{i},1-\sigma_{i}}\mathbb{U}_{\alpha,\sigma_{i}}$.
Furthermore, for any fixed $0<\sigma\leq1$, the density of
$M_{1}(F_{Y_{\sigma}\mathbb{U}_{\alpha,0}})$ is given by, for
$0<y<1,$
\[
\frac{{(\alpha+1)}^{\sigma/\alpha}y^{\sigma-1}\sin(\curpi\sigma
[1-F_{\mathbb{U}_{\alpha,0}}(y)])}
{{[y^{2\alpha+2}+2y^{\alpha+1} {(1-y)}^{\alpha+1}\cos(\curpi
\alpha)+{(1-y)}^{2\alpha+2}]}^{{\sigma}/{(2\alpha)}}}.
\]
\item[(ii)]The finite-dimensional
distribution of
$(\Upsilon^{\ddagger}_{\alpha}(C_{1}),\ldots,\Upsilon^{\ddagger}_{\alpha
}(C_{k}))$
is such that each $\Upsilon^{\ddagger}_{\alpha}(C_{i})$ is
independent and has distribution
\[
\Upsilon^{\ddagger}_{\alpha}(C_{i})\stackrel{d}=G_{1}M_{1}(F_{Y_{\sigma
_{i}}Z^{1/(\alpha+1)}_{\alpha}}) ,
\]
where
$M_{1}(F_{Y_{\sigma_{i}}Z^{1/(\alpha+1)}_{\alpha}})\stackrel{d}=\beta
_{\sigma_{i},1-\sigma_{i}}M_{\sigma_{i}}(F_{Z^{1/(\alpha+1)}_{\alpha}})$.
Furthermore, for any fixed $0<\sigma\leq1$, the density of
$M_{1}(F_{Y_{\sigma}Z^{1/(\alpha+1)}_{\alpha}})$ is given by, for
$x>0,$
\[
\frac{x^{\sigma-1}{(x+1)}^{{\sigma(1+\alpha)}/{\alpha}-1}\sin(\curpi
\sigma[1-F_{Z_{\alpha}}(x^{\alpha+1})])}
{{[x^{2\alpha+2}+2x^{\alpha+1} \cos(\curpi
\alpha)+1]}^{{\sigma}/{(2\alpha)}}} .
\]
\end{enumerate}
\end{thm}
\begin{pf} Statement (i) follows from Theorem 2.2 and Proposition
2.1 in
combination\break with~(\ref{JLP}). Noting the relationship between
$Z^{1/(\alpha+1)}_{\alpha}$ and $\mathbb{U}_{\alpha,0},$ statement
(ii) follows from Theorem~3.1(ii).
\end{pf}

From this, combined with an application of (\ref{sinid}), we
obtain a description for the densities of
$\Upsilon^{\ddagger}_{\alpha}(1)$ and $\Upsilon_{\alpha}(1).$ In
addition, for $\alpha\leq1/2,$ we obtain a description of the
density of $\Upsilon_{\alpha}(2)$ using~(\ref{idd}).
\begin{prop} Let $\Upsilon_{\alpha}(1)$ and
$\Upsilon^{\ddagger}_{\alpha}(1)$ denote GGC random variables with
parameters $(\alpha, F_{\mathbb{U}_{\alpha,0}})$ and
$(\alpha,F_{Z^{1/(\alpha+1}_{\alpha}})$, respectively. Then:
\begin{enumerate}[(iii)]
\item[(i)]$\Upsilon_{\alpha}(1)\stackrel{d}=G_{1}M_{1}(F_{Y_{\alpha
}\mathbb{U}_{\alpha,0}})$,
where
$M_{1}(F_{Y_{\alpha}\mathbb{U}_{\alpha,0}})\stackrel{d}=\beta_{\alpha
,1-\alpha}\mathbb{U}_{\alpha,\alpha}$
has density, for $0<y<1,$
\[
\frac{\sin(\curpi\alpha)}{\curpi}\frac{(\alpha+1)y^{\alpha-1}{(1-y)}^{\alpha
+1}}{{[y^{2\alpha+2}+2y^{\alpha+1}
{(1-y)}^{\alpha+1}\cos(\curpi\alpha)+{(1-y)}^{2\alpha+2}]}} .
\]
\item[(ii)]$\Upsilon^{\ddagger}_{\alpha}(1)\stackrel
{d}=G_{1}M_{1}(F_{Y_{\alpha}Z^{1/(\alpha+1)}_{\alpha}})$,
where $M_{1}(F_{Y_{\alpha}Z^{1/(\alpha+1)}_{\alpha}})$ has
density
\[
\frac{\sin(\curpi\alpha)}{\curpi}\frac{x^{\alpha-1}{(1+x)}^{\alpha
}}{{[x^{2\alpha+2}+2x^{\alpha+1}
\cos(\curpi\alpha)+1]}} \qquad{\mbox{for }}x>0.
\]
\item[(iii)]Supposing that $\alpha\leq1/2,$ then the
GGC$(2\alpha,F_{Z^{1/(\alpha+1}_{\alpha}})$ random variable
$\Upsilon^{\ddagger}_{2\alpha}(1)\stackrel{d}=\Upsilon^{\ddagger}_{\alpha}(2)$
is equivalent in distribution to
$G_{1}M_{1}(F_{Y_{2\alpha}Z^{1/(\alpha+1)}_{\alpha}}),$ where
$M_{1}(F_{Y_{2\alpha}Z^{1/(\alpha+1)}_{\alpha}})$ has density
\[
\frac{2x^{2\alpha-1}{(x+1)}^{2\alpha+1}\sin(\curpi\alpha)[\cos(\curpi\alpha
)+x^{\alpha+1}]}
{{[x^{2\alpha+2}+2x^{\alpha+1} \cos(\curpi\alpha)+1]}^{2}}\qquad{\mbox{for }}x>0.
\]
\end{enumerate}
\end{prop}
%
\subsection{An example connected to Diaconis and Kemperman}
Note that we have the following convergence in distribution
results, as $\alpha\rightarrow0:$
\[
\tilde{M}_{\alpha,\theta}(F_{U})=M_{\theta}(F_{\mathbb{U}_{\alpha
,0}})\stackrel{\mathrm{d}}\rightarrow
M_{\theta}(F_{U})\qquad{\mbox{for }}\theta>0
\]
and
\[
\mathbb{U}_{\alpha,0}\stackrel{\mathrm{d}}\rightarrow1-U.
\]
Furthermore, setting $W=(1-U)/U=G_{1}/G'_{1}$, we have
\[
M_{\theta}(F_{Z^{1/(\alpha+1)}_{\alpha}})\stackrel{\mathrm{d}}\rightarrow
M_{\theta}(F_{W})\quad {\mbox{and}}\quad Z^{1/(\alpha+1)}_{\alpha}\stackrel{\mathrm{d}}\rightarrow W,
\]
where the last statement can be read from Chaumont and
Yor~\cite{Chaumont}, page 155 and page 169. It is then natural to
investigate the laws of the random processes connected with the
GGC$(\theta,F_{U})$ and GGC$(\theta,F_{W})$ random variables. It
is known from Diaconis and Kemperman~\cite{Diaconis} that the
density of $M_{1}(F_{U})$ is
%
\begin{equation} \label{DPUni} \frac{\mathrm{e}}{\curpi}\sin(\curpi
y)y^{-y}{(1-y)}^{-(1-y)}\qquad{\mbox{for }}0<y<1.
\end{equation}
Note,
furthermore, that $\tilde{T}_{1}\stackrel{\mathrm{d}}=G_{1}M_{1}(F_{U})$ is
GGC$(1,F_{U})$ and has Laplace transform
\[
\mathbb{E}\bigl[\mathrm{e}^{-\lambda G_{1}M_{1}(F_{U})}\bigr]=\mathrm{e}^{-\psi_{F_{U}}(\lambda)}=\mathrm{e}{(1+\lambda)}^{-({(\lambda+1)}/{\lambda})}.
\]
Now, $G_{1}M_{1}(F_{W})$ is a GGC$(1,F_{W})$ with
$\psi_{F_{W}}(\lambda)=\frac{\lambda}{\lambda-1}\log(\lambda).$
Theorem 3.1 shows that $M_{1}(F_{U})$ arises from tilting the
density of $G_{1}M_{1}(F_{W}).$ The density of $M_{1}(F_{W})$ is
obtained by applying statement (ii) of Theorem 3.1 to
(\ref{DPUni}), or by statement (ii) of Proposition 3.1, and is
given by
\[
\xi_{F_{W}}(x)=\frac{1}{\curpi}\sin\biggl(\frac{\curpi
x}{1+x}\biggr)x^{-{x}/{(1+x)}}\qquad{\mbox{for }}x>0.
\]

We now apply Theorem 2.2 and Proposition 2.1 to give a description
of the finite-dimensional distribution of the subordinators
associated with the two random variables above.
\begin{prop}Let U denote a uniform [0,1] random variable and let
$W=G_{1}/G'_{1}$
denote a~ratio of independent exponential $(1)$ random variables.
\begin{enumerate}[(iii)]
\item[(i)] Suppose that
$(\tilde{\zeta}_{1}(t); 0<t<1)$ is a GGC$(1,F_{U})$ subordinator.
Then, for $(C_{1},\ldots,C_{k})$, a~disjoint partition of $(0,1),$
the finite-dimensional distribution
has joint
density as in~(\ref{fidi1}), with
\[
g_{\sigma_{i},F_{U}}(x_{i})=\int_{0}^{1}\mathrm{e}^{-{x_{i}}/{y}}\frac{\mathrm{e}^{\sigma_{i}}}{\curpi}\sin\bigl(\curpi
\sigma_{i}(1-
y)\bigr)y^{\sigma_{i}(1-y)-2}{(1-y)}^{-\sigma_{i}(1-y)}\,\mathrm{d}y
\]
for $i=1,\ldots,k.$
\item[(ii)] That is,
$\tilde{\zeta}_{1}(C_{i})\stackrel{d}=G_{1}M_{1}(F_{UY_{\sigma_{i}}})$
and they are independent for $i=1,\ldots, k.$ Furthermore, the density
of
$M_{1}(F_{UY_{\sigma_{i}}})\stackrel{d}=\beta_{\sigma_{i},1-\sigma
_{i}}M_{\sigma_{i}}(F_{U})$
is
\[
\frac{\mathrm{e}^{\sigma_{i}}}{\curpi}\sin\bigl(\curpi\sigma_{i}(1-
y)\bigr)y^{\sigma_{i}(1-y)-1}{(1-y)}^{-\sigma_{i}(1-y)}
\]
for $0<y<1.$
\item[(iii)]If $(\zeta_{1}(t); 0<t<1)$ is a GGC$(1,F_{W})$
subordinator, then
the finite-dimensional distribution
$(\zeta_{1}(C_{1}),\ldots,\zeta_{1}(C_{k}))$ is now described,
with
\[
g_{\sigma_{i},F_{W}}(x_{i})=\int_{0}^{\infty}\mathrm{e}^{-{x_{i}}/{w}}\frac{1}{\curpi}\sin\biggl(\frac{\curpi
\sigma_{i}}{1+w}\biggr)w^{{\sigma_{i}}/{(1+w)}-2}\,\mathrm{d}w.
\]
\item[(iv)] That is,
$\zeta_{1}(C_{i})\stackrel{d}=G_{1}M_{1}(F_{WY_{\sigma_{i}}})$ and
they are independent for $i=1,\ldots, k.$ Furthermore, the density of
$M_{1}(F_{WY_{\sigma_{i}}})\stackrel{d}=\beta_{\sigma_{i},1-\sigma
_{i}}M_{\sigma_{i}}(F_{W})$
is
\[
\frac{1}{\curpi}\sin\biggl(\frac{\curpi
\sigma_{i}}{1+x}\biggr)x^{{\sigma_{i}}/{(1+x)}-1}
\]
for $x>0.$
\end{enumerate}
\end{prop}
\begin{pf} This now follows from Theorem 2.2, Proposition 2.1 and
(\ref{DPUni}). Specifically, note that
$\mathcal{C}(F_{U})=(0,\infty),$ so for any $0<\sigma<1,$
\[
\sin(\curpi F_{UY_{\sigma}}(u))=\sin\bigl(\curpi\sigma(1-u)\bigr)
\]
for $0<u<1$ and $0$ otherwise. Furthermore, from~(\ref{DPUni}), or
by direct argument, it is easy to see that
\[
\Phi_{F_{U}}(y)=-\log\bigl(y^{-y}{(1-y)}^{-(1-y)}\bigr)-1.
\]
This fact is also evident from Diaconis and
Kemperman~\cite{Diaconis}. It follows that
$M_{1}(F_{UY_{\sigma}})$ has density
\[
\frac{\mathrm{e}^{\sigma}}{\curpi}\sin\bigl(\curpi\sigma(1-
y)\bigr)y^{\sigma(1-y)-1}{(1-y)}^{-\sigma(1-y)}\qquad{\mbox{for }}0<y<1.
\]
The density for $M_{1}(F_{WY_{\sigma}})$ is obtained in a similar fashion
by Proposition 3.1.
\end{pf}

\begin{rem}Setting
\[
A_{c}\stackrel{\mathrm{d}}=\frac{cG_{1}}{cG_{1}+G'_{1}} ,
\]
one can easily obtain the density of the random variable
$M_{1}(F_{A_{c}})$ for each $c>0$ by using statement (ii) of Theorem
3.1. Note, also, that one can deduce from the density of
$M_{1}(F_{W})$ that $\Phi_{F_{W}}(x)=[x/(1+x)]\log(x).$ Hence, in
this case, an application of Proposition 3.2 shows that
\[
\Phi_{F_{A_{c}}}(y)=
\frac{y}{c(1-y)+y}\log\biggl(\frac{y}{c(1-y)}\biggr)-\frac{c\log(c)}{c-1}+\log
\bigl(c(1-y)\bigr) .
\]
We note that, otherwise, it is not easy to calculate $\Phi_{A_{c}}$,
in this case, by direct arguments.
\end{rem}

\subsection{The finite-dimensional distribution of subordinators of
Bertoin \textit{et al.}} Our final example shows how one can apply the
results in
Sections 2 and 3 to obtain new results for subordinators
recently studied by Bertoin \textit{et al.}~\cite{BFRY}.
In particular, they investigate properties of the random
variables corresponding to the lengths of excursions of Bessel
processes straddling an independent exponential time, which can be
expressed as
\[
d_{\mathbf{e}}^{(\alpha)} - g_{\mathbf{e}}^{(\alpha)},
\]
where, for any $t>0$,
%
\begin{equation} \label{III17}
g_{t}^{(\alpha)}:= \operatorname{sup} \{s \le t,  R_{s} =0\},\qquad
d_{t}^{(\alpha)} := \operatorname{inf}\{s \ge t,    R_{s}=0\}
\end{equation}
for $(R_{t},   t \ge0)$ a Bessel process starting from $0$ with
dimension $d=2(1- \alpha)$, with $0< d <2$ or, equivalently, $0<
\alpha<1$. Additionally, $\mathbf{e}\stackrel{\mathrm{d}}=G_{1},$ an
exponentially distributed random variable with mean $1.$ See also
Fujita and Yor~\cite{FY2} for closely related work.

In order to avoid confusion, we will now denote relevant random
variables appearing originally as $\Delta_{\alpha}$ and $G_{\alpha}$
in~Bertoin \textit{et al.}~\cite{BFRY} as
$\Sigma_{\alpha}$ and $\mathbb{G}_{\alpha}$, respectively.
From~Bertoin \textit{et al.}~\cite{BFRY}, let
$(\Sigma_{\alpha}(t);t>0)$ denote a subordinator such that
\begin{eqnarray*}
\mathbb{E}\bigl[\mathrm{e}^{-\lambda
\Sigma_{\alpha}(t)}\bigr]&=&{\bigl((\lambda+1)^{\alpha}-\lambda^{\alpha}
\bigr)}^{t}\\
&=&\exp\bigl(-t(1-\alpha)\mathbb{E}[\log(1+\lambda/\mathbb{G}_{\alpha})]\bigr) ,
\end{eqnarray*}
where, from Bertoin \textit{et al.} \cite{BFRY}, Theorems 1.1
and 1.3, $\mathbb{G}_{\alpha}$ denotes a random variable
such that
\[
\mathbb{G}_{\alpha}\stackrel{\mathrm{d}}=\frac{Z^{1/\alpha}_{1-\alpha
}}{1+Z^{1/\alpha}_{1-\alpha}}
\]
and has density on $(0,1)$ given by
\[
f_{\mathbb{G}_{\alpha}}(u)=\frac{\alpha
\sin(\curpi\alpha)}{(1-\alpha)\curpi}\frac{u^{\alpha-1}{(1-u)}^{\alpha
-1}}{u^{2\alpha}-2(1-u)^{\alpha}u^{\alpha}\cos(\curpi\alpha
)+{(1-u)}^{2\alpha}}.
\]
Hence, it follows that the random variable $1/\mathbb{G}_{\alpha}$
takes its values on $(1,\infty)$ with probability one and has cdf
satisfying
\[
1-F_{1/\mathbb{G}_{\alpha}}(x)=F_{\mathbb{G}_{\alpha
}}(1/x)=F_{Z_{1-\alpha}}\bigl({(x-1)}^{-\alpha}\bigr).
\]
As noted by Bertoin \textit{et al.} \cite{BFRY},
$(\Sigma_{\alpha}(t);t>0)$ is a GGC$(1-\alpha,
F_{1/\mathbb{G}_{\alpha}})$ subordinator, where the
GGC$(1-\alpha,F_{1/\mathbb{G}_{\alpha}})$ random variable
$\Sigma_{\alpha}\stackrel{\mathrm{d}}=\Sigma_{\alpha}(1)$ satisfies
\[
\Sigma_{\alpha}\stackrel{\mathrm{d}}=d_{\mathbf{e}}^{(\alpha)} -
g_{\mathbf{e}}^{(\alpha)}\stackrel{\mathrm{d}}=\frac{G_{1-\alpha}}{\beta_{\alpha
,1}}\stackrel{\mathrm{d}}=\frac{G_{1-\alpha}}{U^{1/\alpha}},
\]
where $U$ denotes a uniform$[0,1]$ random variable and, for clarity,
$G_{1-\alpha}$ is a gamma$(1-\alpha,1)$ random variable. This
means that the density of $\Sigma_{\alpha}$ is
\[
\frac{\alpha}{\Gamma(1-\alpha)}x^{-\alpha-1}(1-\mathrm{e}^{-x})\qquad{\mbox{for }}x>0.
\]

It is evident, as investigated in Fujita and Yor~\cite{FY}, that
\[
M_{1-\alpha}(F_{1/\mathbb{G}_{\alpha}})\stackrel{\mathrm{d}}=\frac{1}{\beta
_{\alpha,1}}\stackrel{\mathrm{d}}=U^{-1/\alpha} .
\]
\begin{rem} Note that when $\alpha=1/2,$
$\mathbb{G}_{1/2}\stackrel{\mathrm{d}}=\beta_{1/2,1/2}$. It is known that
for each fixed $t,$
\[
\Sigma_{1/2}(t)\stackrel{\mathrm{d}}=\frac{G_{t/2}}{\beta_{1/2,(1+t)/2}},
\]
which implies that
%
\begin{equation}
M_{t/2}(F_{1/\mathbb{G}_{1/2}})=M_{t/2}(F_{1/\beta_{1/2,1/2}})\stackrel
{\mathrm{d}}=\frac{1}{\beta_{1/2,(1+t)/2}} .
\label{JYbeta}
\end{equation}
This result may be found in James and Yor~\cite{JamesYor}. Related
to this fact, Cifarelli and Melilli~\cite{CifarelliMelilli} have
shown that
$M_{t/2}(F_{\beta_{1/2,1/2}})\stackrel{\mathrm{d}}=\beta_{(t+1)/2,(t+1)/2}$
for $t>0.$
\end{rem}

 In regards to exponentially tilting
GGC$(1-\alpha,F_{1/\mathbb{G}_{\alpha}})$, note that for $c>0,$
\[
\frac{c/\mathbb{G}_{\alpha}}{c/\mathbb{G}_{\alpha}+1}=\frac{c}{\mathbb
{G}_{\alpha}+c}.
\]
Thus, a GGC$(1-\alpha, F_{c/(\mathbb{G}_{\alpha}+c)})$
subordinator, say $(\Sigma^{\dagger}_{\alpha,c}(t),t\leq
1/(1-\alpha))$, arises from exponential tilting. Naturally, the
density of $\Sigma^{\dagger}_{\alpha,c}(1)/c$ is given by
\[
\frac{\alpha x^{-\alpha-1}\mathrm{e}^{-cx}(1-\mathrm{e}^{-x})}{[(c+1)^{\alpha}-c^{\alpha}]\Gamma(1-\alpha)}\qquad{\mbox{for
}}x>0.
\]
Equivalently,
$\Sigma^{\dagger}_{\alpha,c}(1)\stackrel{\mathrm{d}}=G_{1-\alpha}M_{1-\alpha
}(F_{c/(\mathbb{G}_{\alpha}+c)}) ,$
where $M_{1-\alpha}(F_{c/(\mathbb{G}_{\alpha}+c)})$ has density
\[
\frac{\alpha
c^{\alpha}}{(c+1)^{\alpha}-c^{\alpha}}u^{-\alpha-1}\qquad{\mbox{for
}}\frac{c}{c+1}<u<1.
\]

Now, using the facts discussed above, we will show how to use the
results in Section 2 to explicitly describe the finite-dimensional
distribution of the subordinators $(\Sigma_{\alpha}(t), t>0)$ and
$(\Sigma^{\dagger}_{\alpha,c}(t), t>0)$ over the range $0<t\leq
1/(1-\alpha)$.
Additionally, the analysis will also yield expressions for mean
functionals based on $F_{1/\mathbb{G}_{\alpha}}.$ First, note
that, using~(\ref{sinp}), one has
%
\begin{equation}
\label{sinap} \sin(\curpi
F_{Y_{1-\alpha}/\mathbb{G}_{\alpha}}(x))=\cases{
\sin\bigl(\curpi(1-\alpha) F_{\mathbb{G}_{\alpha}}(1/x)\bigr), &\quad  ${\mbox{if }} x>
1,$\vspace*{2pt}\cr
\sin\bigl(\curpi(1-\alpha)\bigr),&\quad ${\mbox{if }} 0<x\leq1,$}
\end{equation}
where, again using the properties of $F_{Z_{1-\alpha}},$ as deduced
from James~\cite{JamesLinnik}, Proposition 2.1(iii),
%
\begin{equation}
\label{cdfid}
\sin\bigl(\curpi(1-\alpha)
F_{\mathbb{G}_{\alpha}}(1/x)\bigr)=\frac{\sin(\curpi(1-\alpha
))}{{[{(x-1)}^{2\alpha}-2{(x-1)}^{\alpha}\cos(\curpi
\alpha)+1]}^{1/2}}.
\end{equation}

We now use this to calculate
%
\begin{equation}
\Phi_{F_{1/\mathbb{G}_{\alpha}}}(x)=\mathbb{E}[\log(|x-1/\mathbb
{G}_{\alpha}|)\indic(x\neq
1/\mathbb{G}_{\alpha})]. \label{phifunction}
\end{equation}
\begin{prop}For $0<\alpha<1,$ consider
$\Phi_{F_{1/\mathbb{G}_{\alpha}}}(x)$ as defined in
(\ref{phifunction}). Then,
%
\begin{eqnarray}
\label{rfunction}
&&\Phi_{F_{1/\mathbb{G}_{\alpha}}}(x)
\nonumber
\\[-8pt]
\\[-8pt]
\nonumber
&&\quad=\cases{
\displaystyle\frac{1}{2(1-\alpha)}\biggl[\log\biggl(\frac{x^{2}}{[(x-1)^{2\alpha}-2(x-1)^{\alpha
}\cos(\curpi\alpha)+1]}\biggr)\biggr], &\quad  ${\mbox{if }} x> 1,$\cr
\displaystyle\frac{1}{1-\alpha}\log\bigl({x}/{[1-{(1-x)}^{\alpha}]}\bigr),&\quad  ${\mbox{if }}
0<x\leq1.$}\hspace*{24pt}
\end{eqnarray}
\end{prop}
\begin{pf}
Using simple beta--gamma algebra, we have
\[
\Sigma_{\alpha}\stackrel{\mathrm{d}}=\frac{G_{1-\alpha}}{\beta_{\alpha,1}}\stackrel
{\mathrm{d}}=G_{1}\frac{\beta_{1-\alpha,\alpha}}{U^{1/\alpha}} .
\]
Hence, applying Theorem 2.1, with $\theta=1$ and
$\sigma=1-\alpha,$ it follows that $\Sigma_{\alpha}$ is also
GGC$(1,F_{Y_{1-\alpha}/\mathbb{G}_{\alpha}})$ and
%
\begin{equation}
B_{\alpha}:=\frac{\beta_{1-\alpha,\alpha}}{\beta_{\alpha,1}}\stackrel
{\mathrm{d}}=\frac{\beta_{1-\alpha,\alpha}}{U^{1/\alpha}}\stackrel
{\mathrm{d}}=M_{1}(F_{Y_{1-\alpha}/\mathbb{G}_{\alpha}}).
\label{key5}
\end{equation}

By standard calculations, the density of
$B_{\alpha}=\beta_{1-\alpha,\alpha}/\beta_{\alpha,1}$ is given by
\[
f_{B_{\alpha}}(x)=\frac{\sin(\curpi
(1-\alpha))}{\curpi}x^{-\alpha-1}[1-{(1-x)}^{\alpha}\indic(x\leq1)] .
\]
However, we see from~(\ref{key5}) that $B_{\alpha}\stackrel{\mathrm{d}}=
M_{1}(F_{Y_{1-\alpha}/\mathbb{G}_{\alpha}}).$ Hence, Theorem 2.2
applies and the density of $B_{\alpha}$ can be written as
\[
f_{B_{\alpha}}(x)=\frac{x^{-\alpha}}{\curpi}\sin(\curpi
F_{Y_{1-\alpha}/\mathbb{G}_{\alpha}}(x))\mathrm{e}^{-(1-\alpha)\Phi_{F_{1/\mathbb{G}_{\alpha}}}(x)} .
\]
Now, equating the two forms of the density of $B_{\alpha}$ and
using (\ref{sinap}) and~(\ref{cdfid}), one then obtains the
expression for $\Phi_{F_{1/\mathbb{G}_{\alpha}}}.$
\end{pf}

Now, for $z>0,$ define the function
\[
\mathcal{S}_{\alpha,\sigma}(z)={\sin(\curpi\sigma
F_{Z_{1-\alpha}}(z^{-\alpha}))}{[z^{2\alpha}-2z^{\alpha}\cos(\curpi\alpha
)+1]}^{{\sigma}/{(2(1-\alpha))}}
\]
and define,
\[
\mathcal{D}_{\alpha,\sigma}(x)=\cases{
\sin(\curpi\sigma){[1-{(1-x)}^{\alpha}]}^
{{\sigma}/{(1-\alpha)}}, &\quad  ${\mbox{if }} x\leq1,$\cr
\mathcal{S}_{\alpha,\sigma}(x-1),& \quad ${\mbox{if }} x> 1.$}
\]

Hereafter, $(C_{1},\ldots,C_{k})$ will denote an arbitrary
disjoint partition of the interval $(0,1/(1-\alpha)]$ with lengths
$|C_{i}|,$ and $\sigma_{i}=(1-\alpha)|C_{i}|$ for $i=1,\ldots,k.$

\begin{thm}Consider the GGC$(1-\alpha,F_{1/\mathbb{G}_{\alpha}})$ subordinator
$(\Sigma_{\alpha}(t),t\leq1/(1-\alpha))$ and, for each fixed
$c>0,$ the GGC$(1-\alpha,F_{c/(\mathbb{G}_{\alpha}+c)}$
subordinator $(\Sigma^{\dagger}_{\alpha,c}(t),t\leq
1/(1-\alpha))$. The following results then hold:
\begin{enumerate}[(ii)]
\item[(i)] The finite-dimensional
distribution of
$(\Sigma_{\alpha}(C_{1}),\ldots,\Sigma_{\alpha}(C_{k}))$ is such
that each $\Sigma_{\alpha}(C_{i})$ is independent and has
distribution
\[
\Sigma_{\alpha}(C_{i})\stackrel{d}=G_{1}M_{1}(F_{Y_{\sigma_{i}}/\mathbb
{G}_{\alpha}}),
\]
where
$M_{1}(F_{Y_{\sigma_{i}}/\mathbb{G}_{\alpha}})\stackrel{d}=\beta_{\sigma
_{i},1-\sigma_{i}}M_{\sigma_{i}}(F_{1/\mathbb{G}_{\alpha}})$.
Furthermore, for any fixed $0<\sigma\leq1$, the density of
$M_{1}(F_{Y_{\sigma}/\mathbb{G}_{\alpha}})$ is given by
\[
\frac{1}{\curpi
}x^{-({\sigma\alpha}/{(1-\alpha)}+1)}\mathcal{D}_{\alpha,\sigma
}(x)\qquad{\mbox{for }}x>0.
\]

\item[(ii)]For the GGC$(1-\alpha,F_{c/(\mathbb{G}_{\alpha}+c)})$
process, $\Sigma^{\dagger}_{\alpha,c}$, it follows that each
\[
\Sigma^{\dagger}_{\alpha,c}(C_{i})\stackrel{d}=G_{1}M_{1}\bigl(F_{Y_{\sigma
_{i}}c/(\mathbb{G}_{\alpha+c})}\bigr),
\]
where for each $0<\sigma\leq1,$
$M_{1}(F_{Y_{\sigma}c/(\mathbb{G}_{\alpha+c})})$ has density
\[
\frac{{[c(1-y)]}^{{\sigma\alpha}/{(1-\alpha)}}
\mathcal{D}_{\alpha}({y}/{(c(1-y))})}{\curpi
{[(c+1)^{\alpha}-c^{\alpha}]}^{\sigma}y^{{\sigma\alpha}/{(1-\alpha)}+1}}
\qquad{\mbox{for }}0<y<1.
\]
\end{enumerate}
\end{thm}
\begin{pf} From Theorem 2.2, we have that the general form of the
density of $M_{1}(F_{Y_{\sigma}/\mathbb{G}_{\alpha}})$ is given by
\[
\frac{x^{\sigma-1}}{\curpi}\sin(\curpi
F_{Y_{\sigma}/\mathbb{G}_{\alpha}}(x))\mathrm{e}^{-\sigma
\Phi_{F_{1/\mathbb{G}_{\alpha}}}(x)} .
\]
The proof is completed by applying Proposition 4.3 and
(\ref{sinap}) and~(\ref{cdfid}).
\end{pf}
\begin{rem}
The process $\Sigma_{\alpha,c}(t)/c$ is well defined for $c\ge0$
and $0\leq\alpha<1,$ and presents itself as an interesting class
worthy of further investigation. Letting $c\rightarrow0$, it is
evident that $\Sigma^{\dagger}_{\alpha,c}(1)/c$ converges to
$\Sigma_{\alpha}(1).$ As shown by Bertoin \textit{et al.} \cite{BFRY}, Section 3.6.3, $\Sigma_{0,c}(1)/c,$ for $c>0$, has
a similar interpretation as $\Sigma_{\alpha}(1),$ but where
the Bessel process $(R_{t},t>0)$ is now replaced by a diffusion
process whose inverse local time at $0$ is distributed as a gamma
subordinator $(\gamma_{l}/c;l>0).$ Furthermore, albeit not
explicitly addressed in Bertoin \textit{et al.} \cite
{BFRY}, the random variable
$\Sigma_{\alpha,c}(1)/c\stackrel{\mathrm{d}}=d_{\mathbf{e}}^{(\alpha,c)} -
g_{\mathbf{e}}^{(\alpha,c)}$ has a similar interpretation where
$(R_{t},t>0)$ is now replaced by a process
$(R^{(\alpha,c)}_{t},t>0)$ whose inverse local time is distributed
as a generalized gamma subordinator, that is, a subordinator
whose L\'{e}vy density is specified by $Cy^{-\alpha-1}\mathrm{e}^{-cy}$ for $y>0.$ This interpretation may be deduced from
Donati-Martin and Yor~(\cite{Donati}, see page 880 (1.c)), where
$R^{(\alpha,c)}$ equates with a downwards Bessel process with
drift $c.$
\end{rem}

Bertoin \textit{et al.}~\cite{BFRY} also show that a
GGC $(1-\alpha,F_{\mathbb{G}_{\alpha}})$ random variable satisfies
\[
G_{1-\alpha}M_{1-\alpha}(F_{\mathbb{G}_{\alpha}})=G_{1-\alpha}U .
\]
Hence, the Laplace transform of the GGC
$(1-\alpha,F_{\mathbb{G}_{\alpha}})$ subordinator, say
$(\mathcal{Z}^{\dagger}_{\alpha, 1}(t),t>0),$ is given by
\[
{\biggl(\frac{1}{\alpha\lambda}[(\lambda+1)^{\alpha}-1]\biggr)}^{t} .
\]
Additionally, using the fact that
%
\begin{equation}
\frac{1}{\mathbb{G}_{\alpha}}\stackrel{\mathrm{d}}=\frac{1}{Z^{1/\alpha}_{1-\alpha
}}+1\stackrel{\mathrm{d}}=Z^{1/\alpha}_{1-\alpha}+1
\label{id3}
\end{equation}
leads to
\[
M_{1-\alpha}(F_{Z^{1/\alpha}_{1-\alpha}})\stackrel{\mathrm{d}}=M_{1-\alpha
}(F_{1/\mathbb{G}_{\alpha}})-1\stackrel{\mathrm{d}}=\frac{G_{1}}{G_{\alpha}},
\]
which leads to a description of a
GGC$(1-\alpha,F_{Z^{1/\alpha}_{1-\alpha}})$ subordinator. The
above points may also be found in the survey paper of James,
Roynette and Yor~\cite{JRY}.
\begin{thm}Consider the
GGC$(1-\alpha,F_{Z^{1/\alpha}_{1-\alpha}})$ subordinator
$(\mathcal{Z}_{\alpha}(t),t\leq1/(1-\alpha))$ and the
GGC$(1-\alpha,F_{\mathbb{G}_{\alpha}})$ subordinator
$(\mathcal{Z}^{\dagger}_{\alpha,1}(t),t\leq1/(1-\alpha)).$ The
following results then hold:
\begin{enumerate}[(ii)]
\item[(i)]The finite-dimensional
distribution of
$(\mathcal{Z}_{\alpha}(C_{1}),\ldots,\mathcal{Z}_{\alpha}(C_{k}))$
is such that each $\mathcal{Z}_{\alpha}(C_{i})$ is independent and
is equivalent in distribution to
\[
\mathcal{Z}_{\alpha}(C_{i})\stackrel{d}=G_{1}M_{1}(F_{Y_{\sigma
_{i}}Z^{1/\alpha}_{1-\alpha}}).
\]
Furthermore, for any fixed $0<\sigma\leq1$, the density of
$M_{1}(F_{Y_{\sigma}Z^{1/\alpha}_{1-\alpha}})\stackrel{d}=\beta_{\sigma
,1-\sigma}M_{\sigma}(F_{Z^{1/\alpha}_{1-\alpha}})$
is given by, for $z>0,$
\[
\frac{z^{\sigma-1}}{
\curpi{(1+z)}^{{\sigma}/{(1-\alpha)}}}\mathcal{S}_{\alpha,\sigma}(z) .
\]
\item[(ii)]Similarly, each
$\mathcal{Z}^{\dagger}_{\alpha}(C_{i})\stackrel
{d}=G_{1}M_{1}(F_{Y_{\sigma_{i}}\mathbb{G}_{\alpha}})$
and, for each fixed $0<\sigma\leq1,$
$M_{1}(F_{Y_{\sigma}\mathbb{G}_{\alpha}})$ has density
\[
\frac{\alpha^{{\sigma}/{(1-\alpha)}}}{
\curpi}{y^{\sigma-1}{(1-y)}^{{\sigma\alpha}/{(1-\alpha)}}
}\mathcal{S}_{\alpha,\sigma}\biggl(\frac{y}{1-y}\biggr) .
\]
\end{enumerate}
\end{thm}
\begin{pf} Apply Theorem 2.2 and Theorem 3.1, where, from (\ref{id3}),
\[
\Phi_{F_{Z^{1/\alpha}_{1-\alpha}}}(z)=\Phi_{F_{1/\mathbb{G}_{\alpha}}}(z+1).
\]
\upqed\end{pf}
%
%
\begin{rem}Note that as $\alpha\rightarrow1,$
\[
M_{\theta}(F_{\mathbb{G}_{\alpha}})\stackrel{\mathrm{d}}\rightarrow
M_{\theta}(F_{U}) \quad{\mbox{and}}\quad
M_{\theta}(F_{Z^{1/\alpha}_{1-\alpha}})\stackrel{\mathrm{d}}\rightarrow
M_{\theta}(F_{W}) .
\]
Hence, they have the same limiting behavior, described in Section
4.2, as the random variables in Section 4.1.
\end{rem}

\section*{Acknowlegements}
This research was supported in
part by Grants HIA05/06.BM03, RGC-HKUST 6159/02P, DAG04/05.BM56
and RGC-HKUST 600907 of the HKSAR.

\printhistory

\end{document}